\documentclass[11pt,a4paper]{amsart}

\usepackage[T1]{fontenc}
\usepackage[utf8]{inputenc}
\usepackage{lmodern}
\usepackage{microtype}
\usepackage[a4paper,margin=30mm]{geometry}
\usepackage{amsmath,amssymb,amsthm}
\usepackage{mathtools}
\usepackage{dsfont}
\usepackage{enumitem}
\usepackage{xcolor}
\usepackage{hyperref}
\usepackage{subfigure}
\usepackage{graphicx}

\usepackage{changes}

\usepackage{mathrsfs}

\numberwithin{equation}{section}

\newtheorem{theorem}{Theorem}[section]
\newtheorem{proposition}[theorem]{Proposition}
\newtheorem{lemma}[theorem]{Lemma}
\newtheorem{corollary}[theorem]{Corollary}

\theoremstyle{definition}
\newtheorem{definition}[theorem]{Definition}

\theoremstyle{remark}
\newtheorem{remark}[theorem]{Remark}

\newcommand{\R}{\mathbb{R}}
\newcommand{\ind}{\mathds{1}}

\newcommand{\cD}{\mathscr{D}}
\newcommand{\cA}{\mathcal{A}}

\newcommand{\cT}{\mathscr{T}}
\newcommand{\cK}{\mathcal{K}}

\title[A kinetic model for microblasting]{A kinetic model for microblasting}

\author[M. Hattori]{Masanari Hattori}
\address{M.H.: Department of Aeronautics and Astronautics, Graduate School of Engineering, Kyoto University,
Kyoto daigaku-katsura, Nishikyo-ku, Kyoto 615-8540, Japan}
\email{hattori.masanari.4r@kyoto-u.ac.jp}

\author[F. Salvarani]{Francesco Salvarani}
\address{F.S.: Dipartimento di Matematica ``F.
Casorati'', Universit\`a degli Studi di Pavia, Italy \& De Vinci Higher Education, De Vinci Research Center, Paris, France}
\email{francesco.salvarani@unipv.it}

\author[S. Takata]{Shigeru Takata}
\address{S.T.: Department of Aeronautics and Astronautics, Graduate School of Engineering, Kyoto University,
Kyoto daigaku-katsura, Nishikyo-ku, Kyoto 615-8540, Japan}
\email{takata.shigeru.4a@kyoto-u.ac.jp}

\begin{document}

\begin{abstract}
We propose a kinetic framework for microblasting that combines particle transport, cumulative material weakening and threshold-activated erosion. The evolving geometry is recovered from the damage field and particle absorption is described through space--time characteristics, allowing the model to handle nonsmooth erosion fronts and topology changes without an explicit interface equation. The positivity of these combined mechanisms yields a monotone structure, which is used to construct a global minimal nonnegative mild solution on every finite time interval. A particle-based numerical method is also developed, and numerical experiments are used to explore the qualitative behavior of the model and its ability to represent evolving erosion patterns.
\end{abstract}

\bibliographystyle{abbrv}
\date{}

\maketitle

\medskip

\noindent\textbf{Keywords.}
Kinetic transport; microblasting; cumulative damage; evolving domains;
absorbing boundary; free boundary.

\medskip

\noindent\textbf{MSC 2020.}
35R37, 35Q84, 35R35, 35F16, 65M25.

\section{Introduction}
\label{sec:introduction}

Microblasting is a high-precision surface-treatment technique in which a material is exposed to a directed flux of energetic particles in order to
produce controlled modifications of its surface properties \cite{UZI2021308}.
Its spatial selectivity makes it particularly suitable for fine-scale structuring, selective material removal, surface texturing,
and delicate cleaning.
Applications arise in semiconductor fabrication, biomedical-device preparation, micro-electromechanical systems, precision
optics, microfluidics, and advanced materials processing, where surface topographies must often be controlled at micrometric or submicrometric
scales \cite{wensink2002fabrication}.

In several implementations, the incident particles are directed toward the target through a suitably designed mask.
The mask selectively exposes prescribed regions of the substrate, allowing localized erosion or
patterning while limiting the thermal and mechanical load on the surrounding material.
These features make microblasting particularly useful for small, fragile, or geometrically complex structures \cite{1190517}.

Despite its technological relevance, the mathematical description of microblasting remains incomplete.
Existing approaches include continuum models based on nonlinear diffusion equations
\cite{MR1059002,MR1911147} and atomic-scale cellular descriptions \cite{Tsuda_2010}.
These methods address important aspects of the process, but they are less suited to a setting in which particle transport, material
weakening, absorption, and the geometry of the accessible region must evolve together.

The main difficulty is that erosion is not merely the motion of a prescribed or initially regular boundary.
The particle distribution modifies the material, the accumulated modification changes the region accessible to
future particles, and this geometric change in turn alters the subsequent transport. 
The process is therefore governed by a feedback mechanism between kinetics, damage, and geometry. Moreover, the accessible region may undergo
nucleation, merging, and other topology changes for which a classical interface description is not naturally available.

The strategy adopted in this paper is to formulate the problem through kinetic transport and cumulative damage
rather than through an explicit evolution law for a smooth interface. The particle population is described at the phase-space level,
following the free-molecular, or Knudsen-gas, perspective \cite{sone2007molecular}. The effect of the particles on the
material is represented by a positive cumulative weakening mechanism, while erosion is introduced separately through an activation threshold.

This separation between weakening and erosion is a central feature of the model.
It allows the material to accumulate damage before becoming accessible, accommodates spatially nonlocal interactions, and avoids the
need to identify erosion with a pointwise boundary velocity.
The evolving geometry is recovered from the damage field itself, rather than imposed through an additional interface equation.

A second key element is the use of a space--time formulation based on absorbed characteristics. Instead of tracking a regular moving boundary, the
model determines which terminal portions of particle trajectories remain inside the accessible phase.
This approach continues to be meaningful when the threshold set is nonsmooth or changes topology.
It also provides a natural treatment of absorption without requiring classical normals, curvature, or differentiability of the erosion front.

This strategy offers several advantages over more traditional alternatives.
Diffusion-based models may suppress the directional and finite-speed nature of the incident particle transport, whereas the kinetic description retains
the characteristic structure of particle propagation.
Classical free-boundary models generally require geometric regularity that may be lost during nucleation or merging events.
Fully atomistic models, on the other hand, may resolve microscopic mechanisms in great detail but are less
convenient for a global analytical theory at device-relevant scales.
The present formulation lies between these descriptions: it preserves the
kinetic organization of the particle flux while treating material weakening and geometric evolution at a continuum level.

The same structure also determines the analytical strategy. Positivity of the particle density and of the weakening mechanism induces an order
relation: greater accessibility produces more particle transport, which produces more damage and therefore enlarges the accessible phase.
Starting from the initially accessible region, this order structure leads to a monotone iteration.
The limiting construction yields a global mild solution on every prescribed finite time interval and selects the minimal nonnegative solution.
We underline that the argument does not require smoothness of the evolving threshold set and
does not rely on a contraction principle near the activation threshold.

The proposed framework is therefore designed primarily for robustness under low geometric regularity.

The analytical strategy also suggests a particle-based numerical method. Computational particles can be transported along approximate
characteristics, absorbed according to the currently accessible phase, and used to update a discrete cumulative damage field.
The supercritical region can then be activated iteratively.
This particle formulation mirrors the monotone structure of the analysis and avoids the need to reconstruct a
smooth moving interface at every step.

The numerical simulations described in this article illustrate several features of the proposed kinetic model. 
For a single normally incident particle beam, the computed rate of decrease of the total residual integrity is consistent with the dependence on the incident particle speed incorporated in the kernel function. 
For a substrate containing an initially isolated hole, the eroded region eventually connects to the hole, illustrating a change in the topology of the particle-accessible region. 
For two obliquely incident beams, two initially separate eroded regions develop and eventually merge, illustrating how the evolution of the particle-accessible region gives rise to an indirect interaction between the beams.

The structure of the article is the following. After introducing the kinetic weakening--erosion model and its space--time accessible phase in Section~\ref{sec:model},
the article develops, in Section~\ref{sec:transport}, the required transport and order properties, and then constructs a global minimal mild solution by monotone iteration
in Section~\ref{sec:existence}.
Section~\ref{s:numerical} presents numerical simulations illustrating the behavior of the proposed kinetic model. 
Finally, Section~\ref{sec:discussion} summarizes the results and discusses some remaining open problems.

\section{The kinetic weakening--erosion model}
\label{sec:model}

In this section, we formulate the kinetic weakening--erosion model and establish the topological and measurability
properties required by its mild formulation.
We first define the damage-dependent accessible phase, formalized as the set of space--time points reachable by blasting particles, and study its stability under monotone approximation.
We then introduce a topological backward entrance time along free-transport characteristics, which encodes absorption at the inaccessible phase.
Finally, we use this entrance time to define the absorbing transport operator and the corresponding notion of global mild solution.

\subsection{Preliminary concepts and results} 

\begin{figure}[bt]
    \centering 
\includegraphics[width=0.75\textwidth]{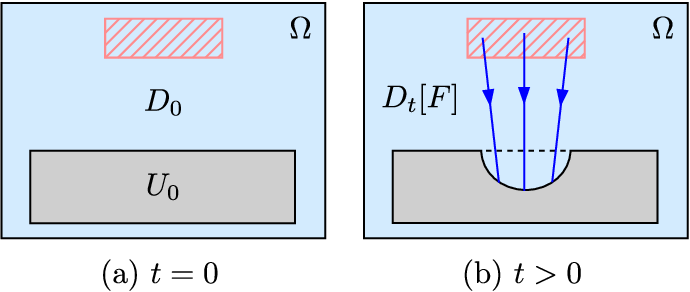}
\caption{
Sketch of the initial and evolved configurations. 
(a) Initial configuration at $t=0$, where $U_0$ is the solid region and $D_0=\Omega\setminus U_0$ is the initially accessible region. 
(b) Configuration at a later time $t>0$, where the accessible region expands to $D_t[F]$ due to erosion; see \eqref{eq:DG}. 
The hatched area indicates the spatial region where particles are generated by the source term $S$. 
In (b), the solid lines with arrowheads represent typical free-flight trajectories of particles, while the dashed line indicates a portion of the boundary of the solid region at the initial time $t=0$.
}
\label{fig:sketch}
\end{figure}

Let $T>0$ and let $\Omega\subset\mathbb R^d$ ($d\geq 1$) be a bounded open set. Let $U_0\subset\Omega$ be a nonempty closed set.
The initially accessible region is $D_0:=\Omega\setminus U_0$ and is assumed to be nonempty  (see Fig.~\ref{fig:sketch}).

Since $\Omega$ is open and $U_0$ is closed in $\mathbb{R}^d$, the set $D_0$ is open in $\mathbb{R}^d$.

The cumulative damage is represented by a nonnegative function
$$
F\, : \, [0,T]\times\overline{\Omega}\to[0,\infty).
$$
We define the residual integrity as
\begin{equation}
\label{eq:uF}
u(t,x):=\ind_{U_0}(x)[1-F(t,x)]_+,
\end{equation}
where $[r]_+:=\max\{r,0\}$.

For $x\in U_0$, the condition $F(t,x)<1$ corresponds to strictly positive residual integrity, whereas $F(t,x)\geq1$
corresponds to vanishing residual integrity.
Outside $U_0$, the integrity is zero by definition.

The state $F=1$ is interpreted as a fully weakened but not yet particle-accessible state. Accessibility occurs after the strict
threshold crossing $F>1$.
This convention is imposed both in the interior and on the initial boundary of the solid.
The strict inequality is essential for the accessible phase to remain open and for it to be stable under increasing pointwise approximation, as shown below.

For a prescribed field
$$
G\, : \, [0,T]\times\overline{\Omega}\to\mathbb{R},
$$
define 
\begin{equation} 
\label{eq:DG}
D_t[G] := D_0\cup\bigl\{x\in U_0:G(t,x)>1\bigr\}, \qquad 0\leq t\leq T.
\end{equation} 
Since $D_0\subset\Omega$ and $U_0\subset\Omega$, one has $D_t[G]\subset\Omega$ for every $t\in[0,T]$.

A first property of $D_t[G]$ is established in the following lemma.

\begin{lemma}
\label{lem:strict-threshold-stability}
Suppose that $G^k(t,x)\to G(t,x)$ pointwise on
$[0,T]\times\overline\Omega$, with
$$
G^k(t,x)\leq G^{k+1}(t,x)
$$
for every $k\in\mathbb N$ and every
$(t,x)\in[0,T]\times\overline\Omega$. Then
$$
D_t[G]=\bigcup_{k\geq0}D_t[G^k]
\qquad\text{for every }t\in[0,T].
$$
\end{lemma}

\begin{proof}
Since $G^k\leq G$, one has
$D_t[G^k]\subseteq D_t[G]$ for every $k$.

Conversely, let $x\in D_t[G]$. If $x\in D_0$, then $x\in D_t[G^k]$ for every $k$. Otherwise, $x\in U_0$ and
$G(t,x)>1$. Since $G^k(t,x)\to G(t,x)$, there exists $k$ such that
$G^k(t,x)>1$. Hence $x\in D_t[G^k]$.
\end{proof}

The particle density is a nonnegative function
$f\, : \, (0,T)\times\Omega\times\mathbb{R}^d\to[0,\infty)$.
Let
$$
K\, : \, \overline\Omega\times\Omega\times\mathbb{R}^d\to[0,\infty)
$$
be a nonnegative measurable kernel. Whenever the following integral is finite, the cumulative damage generated by $f$ is defined by
$$
(\mathcal{K} f)(t,x) := \int_0^t\int_\Omega\int_{\mathbb{R}^d} K(x,x',v)f(s,x',v)\,dv\,dx'\,ds.
$$
The transport-damage model requires
$$
F=\mathcal Kf.
$$

The space--time accessible phase is defined as follows.
For a continuous field $G\in C([0,T]\times\overline\Omega)$, define
$$
\mathscr D[G] :=
\big\{(t,x)\in[0,T]\times\mathbb{R}^d:x\in D_t[G] \big\}.
$$
Equivalently,
\begin{equation}
\label{eq:space-time-accessible-expanded}
\mathscr D[G] = ([0,T]\times D_0)
\cup
\bigl\{(t,x)\in[0,T]\times U_0:G(t,x)>1\bigr\}.
\end{equation}

The interval $[0,T]$ is endowed with its relative topology, and $\mathscr D[G]$ is regarded as a subset of the relative space $[0,T]\times\mathbb{R}^d$.

The following lemma deals with the openness of the accessible phase.
\begin{lemma}
\label{lem:accessible-phase-open}
Let $G\in C([0,T]\times\overline\Omega)$. Then $\mathscr D[G]$ is relatively open in
$[0,T]\times\mathbb R^d$.
In particular, $D_t[G]$ is open in $\mathbb R^d$ for every $t\in[0,T]$.
\end{lemma}

\begin{proof}
Since $D_0=\Omega\setminus U_0$, definition \eqref{eq:DG} gives
$$
\mathbb R^d\setminus D_t[G] =
(\mathbb R^d\setminus\Omega) \cup \bigl\{x\in U_0:G(t,x)\leq1\bigr\}.
$$
The set $\mathbb R^d\setminus\Omega$ is closed. Moreover, $U_0$ is
closed in $\mathbb R^d$, and the continuity of $G(t,\cdot)$ on $\overline\Omega$ implies that
$\bigl\{x\in U_0:G(t,x)\leq1\bigr\}$ is closed in $\mathbb R^d$.
Therefore $\mathbb R^d\setminus D_t[G]$ is closed, and hence $D_t[G]$ is open.

For the space--time phase, its relative complement in $[0,T]\times\mathbb R^d$ is
$$
([0,T]\times\mathbb R^d)\setminus\mathscr D[G]
= [0,T]\times(\mathbb R^d\setminus\Omega)
\cup \bigl\{(t,x)\in[0,T]\times U_0:G(t,x)\leq1\bigr\}.
$$
The first set is relatively closed in $[0,T]\times\mathbb R^d$.
The second one is relatively closed because $[0,T]\times U_0$ is closed in
$[0,T]\times\mathbb R^d$ and $G$ is continuous on $[0,T]\times\overline\Omega$.

Thus the relative complement of $\mathscr D[G]$ is closed, and $\mathscr D[G]$ is relatively open.
\end{proof}

Note that, for a general prescribed continuous field $G$, no temporal monotonicity is assumed, and therefore arbitrary changes of spatial topology may occur.

For the damage field $F=\mathcal Kf$, temporal monotonicity and continuity will be established under the assumptions introduced in Section~\ref{sec:transport}.

The model being based on a transport equation, the absorbed characteristics and entrance time are crucial concepts for our study.

For $(t,x,v)\in[0,T]\times\mathbb{R}^d\times\mathbb{R}^d$, define the backward free-transport characteristic
$$
X_{t,x,v}(s) := x-(t-s)v, \qquad 0\leq s\leq t.
$$

As established in Lemma~\ref{lem:accessible-phase-open}, the space--time accessible phase $\mathscr{D}[G]$ associated with any continuous field $G$ is relatively open. Since the continuity of the damage field $F = \mathcal{K}f$ will be ensured in Section~\ref{sec:Est-ord} (Lemma~\ref{lem:damage-continuity}), we abstract this topological property and define particle transport along characteristics over a general relatively open set $\mathcal{A} \subset [0, T] \times \mathbb{R}^d$.
The characteristic is considered accessible at time $s$ if
$(s,X_{t,x,v}(s))\in\mathcal A$.
In particular, if $\mathcal A\subset[0,T]\times\Omega$, any portion of the characteristic outside $\Omega$ is inaccessible.

\begin{lemma}
\label{lem:terminal-component}
Let $\mathcal A\subset[0,T]\times\mathbb R^d$ be relatively open.
If $(t,x)\in\mathcal A$, then 
$$
I_{\mathcal A}(t,x,v) := \{s\in[0,t]:(s,X_{t,x,v}(s))\in\mathcal A\}
$$
is relatively open in $[0,t]$. Its connected component containing $t$ is either $[0,t]$ or $(\tau,t]$ for a unique $\tau\in[0,t)$. \end{lemma}

\begin{proof}
The map
$$
[0,t]\ni s\mapsto(s,X_{t,x,v}(s))\in[0,T]\times\mathbb{R}^d
$$
is continuous. Since $\mathcal A$ is relatively open, its inverse image $I_{\mathcal A}(t,x,v)$ is relatively open in $[0,t]$.

Every connected component of a relatively open subset of a real interval is an interval that is relatively open in that interval.
Because the component under consideration contains the right endpoint $t$, it has the form $[0,t]$ or $(\tau,t]$ for a unique $\tau\in[0,t)$.
\end{proof}

We now introduce the backward entrance time, which is crucial in the next steps of our study.
\begin{definition}
Denote by $\tau_{\mathcal A}(t,x,v)$ the backward entrance time.
If $(t,x)\notin\mathcal A$, set $\tau_{\mathcal A}(t,x,v):=t$.
If $(t,x)\in\mathcal A$, define $\tau_{\mathcal A}(t,x,v)$ as the left endpoint of the connected component containing $t$ identified
in the preceding lemma, with the convention $\tau_{\mathcal A}=0$ when that component is $[0,t]$.
\end{definition}

Thus, $\tau_{\mathcal A}(t,x,v)$ is the left endpoint of the maximal final time interval over which the characteristic remains in $\mathcal A$.
For $G\in C([0,T]\times\overline\Omega)$, write $\tau_G:=\tau_{\mathscr D[G]}$.

We next prove that backward entrance times are measurable, as deduced in the following lemma.
\begin{lemma}
\label{lem:tau-measurable}
Let $\mathcal A\subset[0,T]\times\mathbb{R}^d$ be relatively open. Then the map
$$
(t,x,v)\mapsto\tau_{\mathcal A}(t,x,v)
$$
is Borel measurable on $[0,T]\times\mathbb{R}^d\times\mathbb{R}^d$. Consequently, $\tau_G$ is Borel measurable whenever $G\in C([0,T]\times\overline\Omega)$.
\end{lemma}

\begin{proof}
We prove that every strict sublevel set of $\tau_{\mathcal A}$ is a Borel set. Fix $a\in(0,T]$. For $q\in\mathbb{Q}\cap[0,a)$, define
$$
E_q := \left\{(t,x,v): q\leq t \ \text{and}\ 
(s,X_{t,x,v}(s))\in\mathcal A \text{ for every } s \in [q,t] \right\}.
$$
We claim that
\begin{equation}
\label{eq:tau-sublevel}
\{\tau_{\mathcal A}<a\} = \{(t,x,v):t<a\} \cup \bigcup_{q\in\mathbb{Q}\cap[0,a)} E_q.
\end{equation}

If $t<a$, then $0\leq\tau_{\mathcal A}(t,x,v)\leq t<a$, regardless of whether $(t,x)$ belongs to $\mathcal A$.

Suppose next that $t\geq a$ and $\tau_{\mathcal A}(t,x,v)<a$. Then $\tau_{\mathcal A}(t,x,v)<t$, so $(t,x)\in\mathcal A$. Choose $q\in\mathbb{Q}$ such that
$$
\tau_{\mathcal A}(t,x,v)<q<a.
$$
By Lemma~\ref{lem:terminal-component}, the characteristic segment over $[q,t]$ lies in $\mathcal A$.
Hence $(t,x,v)\in E_q$.

Conversely, if $(t,x,v)\in E_q$, then the terminal accessible component contains $[q,t]$. Therefore
$$
\tau_{\mathcal A}(t,x,v)\leq q<a.
$$
This proves \eqref{eq:tau-sublevel}.

It remains to prove that $E_q$ is Borel. If $\mathcal A=[0,T]\times\mathbb{R}^d$, then
$$
E_q=\{(t,x,v):t\geq q\},
$$
which is closed. Assume henceforth that the relative complement
$$
\mathcal F:=([0,T]\times\mathbb{R}^d)\setminus\mathcal A
$$
is nonempty. Since $\mathcal A$ is relatively open, $\mathcal F$ is relatively closed in $[0,T]\times\mathbb{R}^d$.
Since $[0,T]\times\mathbb{R}^d$ is closed in $\mathbb{R}^{1+d}$, the set $\mathcal F$ is closed in $\mathbb{R}^{1+d}$.

For $t\geq q$, introduce the fixed-interval parametrization
$$
s=q+r(t-q), \qquad 0\leq r\leq1,
$$
and define
$$
\Phi_q(r,t,x,v) := \operatorname{dist}\!\Bigg( \Big (q+r(t-q),X_{t,x,v}(q+r(t-q))\Big), \mathcal{F} \Bigg).
$$
The distance from a nonempty closed set is continuous, and the map inside the distance is continuous. Hence $\Phi_q$ is continuous on
$[0,1]\times\{(t,x,v):t\geq q\}$.
Define
$$
m_q(t,x,v):=\min_{r\in[0,1]}\Phi_q(r,t,x,v).
$$
Because the minimum is taken over the fixed compact set $[0,1]$, the function $m_q$ is continuous. The segment over $[q,t]$ is contained in the open set $\mathcal A$ if and only if its compact image has strictly positive distance from $\mathcal F$. Therefore
$$
E_q = \{(t,x,v):t\geq q,\ m_q(t,x,v)>0\},
$$
which is a Borel set.

Formula \eqref{eq:tau-sublevel} shows that $\{\tau_{\mathcal A}<a\}$ is Borel for every $a\in(0,T]$. For $a\leq0$ this sublevel set is empty, and for $a>T$ it is the entire domain.
Hence $\tau_{\mathcal A}$ is Borel measurable.
Now let $G\in C([0,T]\times\overline\Omega)$. By Lemma~\ref{lem:accessible-phase-open}, the set $\mathscr D[G]$ is relatively open in $[0,T]\times\mathbb{R}^d$.
Applying the first part of the lemma with $\mathcal A=\mathscr D[G]$, we conclude that $$ \tau_G=\tau_{\mathscr D[G]} $$ is Borel measurable.
\end{proof}

We next prove that  entrance times are monotone and stable.

\begin{lemma}
\label{lem:time-order}
Let $\mathcal A_1,\ \mathcal A_2\subset[0,T]\times\mathbb R^d$ be relatively open and suppose that 
$$
\mathcal A_1\subseteq\mathcal A_2. $$
Then, for every $(t,x,v)\in[0,T]\times\mathbb R^d\times\mathbb R^d$,
$$
\tau_{\mathcal A_2}(t,x,v) \leq \tau_{\mathcal A_1}(t,x,v).
$$
\end{lemma}

\begin{proof}
Every final segment of the characteristic ending at time $t$ that is contained in $\mathcal A_1$ is also contained in $\mathcal A_2$.
Hence, enlarging the accessible set can only enlarge the connected component of accessible times containing $t$ and therefore decrease its left endpoint.
\end{proof}

\begin{lemma}
\label{lem:time-stability}
Let $\{\mathcal A^k\}_{k\geq 0}$ be a sequence of relatively open subsets of $[0,T]\times\mathbb{R}^d$ such that
$$
\mathcal A^k\subseteq\mathcal A^{k+1} \qquad\text{for every }k\geq 0.
$$
Set
$$
\mathcal A:=\bigcup_{k\geq 0}\mathcal A^k, \qquad \tau_k:=\tau_{\mathcal A^k},
\qquad \tau:=\tau_{\mathcal A}.
$$
Then, for every $(t,x,v)$, the sequence $\{\tau_k(t,x,v)\}_{k\geq 0}$ is nonincreasing and
\begin{equation}
\label{eq:tau-monotone-limit}
\lim_{k\to\infty}\tau_k(t,x,v)=\tau(t,x,v).
\end{equation}
\end{lemma}

\begin{proof}
By Lemma
\ref{lem:time-order},
the inclusions
$$
\mathcal A^k\subseteq\mathcal A^{k+1}\subseteq\mathcal A
$$ 
imply
$$
\tau_{k+1}(t,x,v)\leq\tau_k(t,x,v) \quad\text{and}\quad \tau_k(t,x,v)\geq\tau(t,x,v)
$$
for every $k\in\mathbb N$.
Thus the nonincreasing sequence $\{\tau_k(t,x,v)\}$ has a limit $\overline\tau$ satisfying $\overline\tau\geq\tau(t,x,v)$.

If $\tau(t,x,v)=t$, then
$$
t\geq\tau_k(t,x,v)\geq t,
$$
so $\tau_k(t,x,v)=t$ for every $k\in \mathbb N$.

Assume now that $\tau(t,x,v)<t$. Fix $s\in(\tau(t,x,v),t)$. By Lemma~\ref{lem:terminal-component} and the definition of $\tau_{\mathcal A}$, the compact characteristic segment
$$
\mathcal C_s : = \big\{(r,X_{t,x,v}(r)):r\in[s,t]\big\}
$$
is contained in $\mathcal A$.

Since $\mathcal C_s$ is compact and $\mathcal A=\bigcup_{k\geq0}\mathcal A^k$, compactness yields indices $k_1,\ldots,k_m$ such that
$$
\mathcal C_s\subseteq \mathcal A^{k_1}\cup\cdots\cup\mathcal A^{k_m}.
$$
Since the sequence is increasing, setting $N:=\max\{k_1,\ldots,k_m\}$ gives $\mathcal C_s\subseteq\mathcal A^N$.

It follows that $\tau_N(t,x,v)\leq s$. Since the sequence of entrance times is nonincreasing,
$$
\tau_k(t,x,v)\leq s \qquad\text{for every }k\geq N.
$$
Passing to the limit gives $\overline\tau\leq s$. Since this holds for every $s>\tau(t,x,v)$, we obtain $\overline\tau\leq\tau(t,x,v)$. Together with the reverse inequality, this proves \eqref{eq:tau-monotone-limit}.
\end{proof}

\subsection{The absorbing transport operator}

In what follows, whenever a function defined on $\Omega\times\mathbb{R}^d$ is evaluated at a spatial point outside $\Omega$, it is understood to be extended by zero in the spatial variable.
Likewise, time-dependent functions defined on $(0,T)\times\Omega\times\mathbb{R}^d$ are extended by zero in the spatial variable to $(0,T)\times\mathbb{R}^d\times\mathbb{R}^d$. 
Values at the temporal endpoints $s=0$ and $s=T$ may be assigned arbitrarily, since they do not affect Lebesgue integrals in time.

Let
$$
f^{\mathrm{in}}\, : \, \Omega\times\mathbb{R}^d\to[0,\infty) \quad\text{and}\quad S\, : \, (0,T)\times\Omega\times\mathbb{R}^d\to[0,\infty)
$$
be Lebesgue-measurable functions that are finite almost everywhere.

For a relatively open phase $\mathcal A\subset[0,T]\times\mathbb{R}^d$ satisfying $[0,T]\times D_0\subseteq\mathcal A$, define
\begin{align}
\label{eq:TA}
\mathscr T_{\mathcal A}(f^{\mathrm{in}},S)(t,x,v) :={}&
\ind_{\mathcal A}(t,x) \ind_{\{\tau_{\mathcal A}(t,x,v)=0\}}
\ind_{D_0}(x-tv) f^{\mathrm{in}}(x-tv,v) \notag\\
&+ \ind_{\mathcal A}(t,x) \int_{\tau_{\mathcal A}(t,x,v)}^t S\big(s,x-(t-s)v,v\big)\,ds.
\end{align}
For $G\in C([0,T]\times\overline\Omega)$, let
$$
\mathscr T_G:=\mathscr T_{\mathscr D[G]}.
$$
Under the assumption $[0,T]\times D_0\subseteq\mathcal A$, the previous formula implies that a characteristic carries its initial value only if its complete backward trajectory from time $t$ to time $0$ remains accessible.
Indeed, if $x-tv\in D_0$, then $(0,x-tv)\in\mathcal A$, while the condition $\tau_{\mathcal A}(t,x,v)=0$ guarantees accessibility at every positive time along the backward characteristic.
If the final accessible interval along the characteristic begins at a positive time, then the incoming value at that entrance time is zero.

\begin{remark}
\label{rem:no-boundary-regularity}
A key advantage of the present formulation is that no geometric regularity of the moving boundary $\partial D_t[F]$ is required.
The construction uses only the relative openness of the space--time accessible phase $\mathscr D[F]$.
In particular, it does not rely on a parametrization of the evolving interface, on the existence of boundary normals, or on kinetic traces on a classical inflow boundary.
Absorption is encoded intrinsically through the left endpoint of the maximal final segment of each backward characteristic contained in the accessible phase.
The formulation therefore remains meaningful when the accessible region undergoes changes of topology or develops geometrically irregular boundaries.
\end{remark}

The following lemma holds.
\begin{lemma}
\label{lem:characteristic-null-sets}
Let $N_0\subset\mathbb{R}^{2d}$ and
$N_S\subset(0,T)\times\mathbb{R}^{2d}$ be Lebesgue-null sets. Then
$$
\big\{(t,x,v):(x-tv,v)\in N_0\big\}
$$
is null in $(0,T)\times\mathbb{R}^{2d}$, and
$$
\big\{(s,t,x,v):(s,x-(t-s)v,v)\in N_S\big\}
$$
is null in $(0,T)^2\times\mathbb{R}^{2d}$.
\end{lemma}

\begin{proof}
For fixed $t$, the maps $(x,v)\mapsto(x-tv,v)$, and $(s,x,v)\mapsto(s,x-(t-s)v,v)$,
are invertible $C^1$ transformations with Jacobian determinant one.
Their inverse images preserve null sets, and the conclusion follows from Tonelli's theorem applied with respect to $t$.
\end{proof}
We next prove that the transport operator is independent of the choice of measurable representatives, as shown in the following lemma.

\begin{lemma}
\label{lem:representatives}
Let $\mathcal A\subset[0,T]\times\mathbb{R}^d$ be relatively open and satisfy $[0,T]\times D_0\subseteq\mathcal A$. Suppose that
$$
f^{\mathrm{in}}_1=f^{\mathrm{in}}_2 \quad\text{a.e. on }\Omega\times\mathbb{R}^d
$$
and $S_1=S_2$ a.e. on $(0,T)\times\Omega\times\mathbb{R}^d$.
Then
$$
\mathscr T_{\mathcal A}(f^{\mathrm{in}}_1,S_1)
= \mathscr T_{\mathcal A}(f^{\mathrm{in}}_2,S_2)
$$
for almost every $(t,x,v)\in(0,T)\times\mathbb{R}^d\times\mathbb{R}^d$.
\end{lemma}

\begin{proof}
After the prescribed zero extensions, let
$$
N_0:=\{(y,v):f^{\mathrm{in}}_1(y,v)\neq f^{\mathrm{in}}_2(y,v)\},
\qquad N_S:=\{(s,y,v):S_1(s,y,v)\neq S_2(s,y,v)\}.
$$
Both sets are null.
Lemma~\ref{lem:characteristic-null-sets} implies that the initial terms agree for almost every $(t,x,v)$ and that the source
integrands agree for almost every $(s,t,x,v)$.
By Fubini's theorem, for almost every $(t,x,v)$ the latter equality holds for almost every $s$;
hence the corresponding source integrals are equal.
The remaining indicators depend only on $\mathcal A$, $\tau_{\mathcal A}$, and $D_0$.
\end{proof}

We now prove that the transport operator is measurable.
\begin{lemma}
\label{lem:transport-measurable}
Let $\mathcal A\subset[0,T]\times\mathbb{R}^d$ be relatively open and satisfy $[0,T]\times D_0\subseteq\mathcal A$.
If $f^{\mathrm{in}}$ and $S$ are nonnegative Lebesgue-measurable functions,
then $\mathscr T_{\mathcal A}(f^{\mathrm{in}},S)$ admits a Lebesgue-measurable representative.
Consequently, $\mathscr T_G(f^{\mathrm{in}},S)$ admits a Lebesgue-measurable representative for every $G\in C([0,T]\times\overline\Omega)$.
\end{lemma}

\begin{proof}
Choose nonnegative Borel-measurable functions $\widetilde f^{\mathrm{in}}$ and $\widetilde S$ such that
$$
\widetilde f^{\mathrm{in}}=f^{\mathrm{in}}
\quad\text{a.e.}, \qquad
\widetilde S=S \quad\text{a.e.}
$$
after extending to zero when needed.
The maps
$$
(t,x,v)\mapsto(x-tv,v)
\quad
\text{and}
\quad
(s,t,x,v)\mapsto(s,x-(t-s)v,v)
$$
are continuous. Therefore
$$
(t,x,v)\mapsto\widetilde f^{\mathrm{in}}(x-tv,v)
$$
and
$$
(s,t,x,v)\mapsto\widetilde S(s,x-(t-s)v,v)
$$
are Borel measurable.

The indicator $\ind_{\mathcal A}(t,x)$ is Borel measurable. By Lemma~\ref{lem:tau-measurable}, $\tau_{\mathcal A}$ is Borel measurable,
and therefore $\ind_{\{\tau_{\mathcal A}=0\}}$ is Borel measurable. Since $D_0$ is open, the map
$$
(t,x,v)\mapsto\ind_{D_0}(x-tv)
$$
is Borel measurable.

Moreover,
$$
\int_{\tau_{\mathcal A}(t,x,v)}^t
\widetilde S\big(s,x-(t-s)v,v\big)\,ds =
\int_{(0,T)} \ind_{\{\tau_{\mathcal A}(t,x,v)<s<t\}} \widetilde S\big(s,x-(t-s)v,v\big)\,ds.
$$
The integrand on the right-hand side is jointly Borel measurable and nonnegative. Tonelli's theorem implies that its integral with respect to $s$ is measurable in $(t,x,v)$. Hence
$\mathscr T_{\mathcal A}(\widetilde f^{\mathrm{in}},\widetilde S)$ is Borel measurable. By Lemma~\ref{lem:representatives},
$$
\mathscr T_{\mathcal A}(f^{\mathrm{in}},S)
= \mathscr T_{\mathcal A}(\widetilde f^{\mathrm{in}},\widetilde S)
\quad\text{a.e.}
$$
Thus $\mathscr T_{\mathcal A}(f^{\mathrm{in}},S)$ admits a Lebesgue-measurable representative.
\end{proof}
We conclude this section by defining the concept of global mild solution for our problem.
Here and below, the term global refers to existence on the whole prescribed interval $[0,T]$, as opposed to existence only on a smaller time interval.

\begin{definition}
\label{def:solution}
A pair $(F,f)$ is a global mild solution on $[0,T]$ if the following conditions hold:
\begin{enumerate}[label=\textup{(\roman*)}]
\item
$f\in L^1\big((0,T)\times\Omega\times\mathbb{R}^d\big)
\cap L^\infty\big((0,T)\times\Omega\times\mathbb{R}^d\big)$, $f\geq 0$ a.e.;
\item
The cumulative damage field associated with $f$, defined pointwise by
\begin{equation}
\label{eq:canonical-F}
F(t,x) := (\cK f)(t,x):= \int_0^t\int_\Omega\int_{\mathbb{R}^d} K(x,x',v)f(s,x',v)\,dv\,dx'\,ds
\end{equation}
is finite for every $(t,x)\in[0,T]\times\overline\Omega$ and belongs to $C([0,T]\times\overline\Omega)$;
\item
\begin{equation*}
f=\mathscr T_F(f^{\mathrm{in}},S) \qquad\text{a.e. on }(0,T)\times\Omega\times\mathbb{R}^d.
\end{equation*}
\end{enumerate}
The physical integrity and the accessible phase are recovered from \eqref{eq:uF} and \eqref{eq:DG}, respectively.
\end{definition}

\section{Assumptions and transport estimates}
\label{sec:transport}

For the proof of the main result, we impose the following assumptions, which also summarize the geometric setting of Section~\ref{sec:model}.
Unless otherwise stated, Assumptions~\ref{H1}--\ref{H5} remain in force throughout this section.

\begin{enumerate}[label=\textup{(H\arabic*)},ref=\textup{H\arabic*}]
\item
\label{H1}
The set $\Omega\subset\mathbb R^d$ is bounded and open.
The initial solid region is a nonempty closed set $U_0\subset\Omega$.
The initially accessible region is $D_0:=\Omega\setminus U_0$ and is assumed to be nonempty.

\item
\label{H2}
The initial datum satisfies
$$
f^{\mathrm{in}}\in L^1(\Omega\times\R^d)\cap L^\infty(\Omega\times\R^d),
\qquad
f^{\mathrm{in}}\geq0\quad\text{a.e.}
$$
Since the transport operator depends on the initial datum only through
$\ind_{D_0}f^{\mathrm{in}}$, replacing $f^{\mathrm{in}}$ by
$\ind_{D_0}f^{\mathrm{in}}$ entails no loss of generality. We therefore assume that
$$
f^{\mathrm{in}}=\ind_{D_0}f^{\mathrm{in}}
\qquad\text{a.e. on }\Omega\times\R^d.
$$

\item
\label{H3}
The source satisfies $S\geq0$ almost everywhere and
$$
S\in L^1\bigl(0,T;L^\infty(\Omega\times\R^d)\bigr)
\cap L^\infty\bigl(0,T;L^1(\Omega\times\R^d)\bigr).
$$
A nonnegative jointly measurable representative of $S$ is fixed once and for all.

\item
\label{H4}
The kernel is a fixed nonnegative jointly measurable function
$$
K\, : \, \overline\Omega\times\Omega\times\R^d\to[0,\infty)
$$
whose equivalence class with respect to product Lebesgue measure belongs to
$L^\infty(\overline\Omega\times\Omega\times\R^d)$.
Its essential bound is denoted by $\|K\|_{L^\infty}$.

\item
\label{H5}
There exist a nondecreasing function
$$
\omega_K\, : \, [0,\infty)\to[0,\infty),
\qquad
\lim_{r\to0^+}\omega_K(r)=0,
$$
and a measurable set $N_K\subset\Omega\times\R^d$, $|N_K|=0$, such that, for every
$(x',v)\in(\Omega\times\R^d)\setminus N_K$ and every
$x_1,x_2\in\overline\Omega$,
\begin{equation}
\label{eq:kernel-modulus}
|K(x_1,x',v)-K(x_2,x',v)|\leq\omega_K(|x_1-x_2|).
\end{equation}
\end{enumerate}

\begin{lemma}
\label{lem:kernel-representative}
Under Assumptions~\ref{H4}--\ref{H5}, after enlarging $N_K$ by a set of measure zero if necessary, one has
\begin{equation}
\label{eq:kernel-pointwise-bound}
0\leq K(x,x',v)\leq\|K\|_{L^\infty}
\end{equation}
for every $x\in\overline\Omega$ and every
$(x',v)\in(\Omega\times\R^d)\setminus N_K$.
\end{lemma}

\begin{proof}
By Assumption~\ref{H4}, the set
$$
E:=\bigl\{(x,x',v)\in\overline\Omega\times\Omega\times\R^d: K(x,x',v)>\|K\|_{L^\infty}\bigr\}
$$
has product Lebesgue measure zero. By Fubini's theorem, there exists a measurable set
$N_K'\subset\Omega\times\R^d$ of measure zero such that, for every
$(x',v)\notin N_K'$, one has 
$K(x,x',v)\leq\|K\|_{L^\infty}$ for almost every $x\in\Omega$.

Replace $N_K$ by $N_K\cup N_K'$ and fix $(x',v)\notin N_K$. Then
$$
E_{x',v}:=\bigl\{x\in\Omega:K(x,x',v)\leq\|K\|_{L^\infty}\bigr\}
$$
has full measure in $\Omega$ and is therefore dense in $\Omega$. Its closure contains
$\overline\Omega$. Given $x\in\overline\Omega$, choose $x_n\in E_{x',v}$ such that
$x_n\to x$. By Assumption~\ref{H5},
$$
|K(x,x',v)-K(x_n,x',v)|\leq\omega_K(|x-x_n|)\to0.
$$
Consequently,
$$
K(x,x',v)=\lim_{n\to\infty}K(x_n,x',v)\leq\|K\|_{L^\infty}.
$$
The lower bound follows from the pointwise nonnegativity of $K$.
\end{proof}

In what follows, the two assumptions on $S$ will be used separately: the
$L^1_tL^\infty_{x,v}$ bound controls the pointwise size of the transport solution,
whereas the $L^\infty_tL^1_{x,v}$ bound controls its mass.
Unless otherwise specified, the norms $\|\cdot\|_{L^1}$ and $\|\cdot\|_{L^\infty}$ below are taken over $\Omega\times\R^d$.

\subsection{Estimates and order properties\label{sec:Est-ord}}

Set
$$
M_\infty:=\|f^{\mathrm{in}}\|_{L^\infty} +\|S\|_{L^1(0,T;L^\infty)}.
$$

\begin{lemma}
\label{lem:Linf}
Let $\mathcal A\subset[0,T]\times\mathbb R^d$ be relatively open and satisfy
$[0,T]\times D_0\subseteq\mathcal A$. Then
$$
0\leq\cT_{\cA}(f^{\mathrm{in}},S)(t,x,v)\leq M_\infty
$$
for almost every $(t,x,v)\in(0,T)\times\R^d\times\R^d$.
\end{lemma}

\begin{proof}
Positivity follows from \eqref{eq:TA}. By changing the fixed representatives on null sets,
which does not alter the equivalence class of the transport operator, we may assume that
$$
0\leq f^{\mathrm{in}}\leq\|f^{\mathrm{in}}\|_{L^\infty}, \qquad 0\leq S(s,\cdot,\cdot)\leq\|S(s)\|_{L^\infty}
$$
outside suitable null sets. By Lemma~\ref{lem:characteristic-null-sets} and Fubini's theorem,
these bounds remain valid along the characteristic maps for almost every $(t,x,v)$. Hence
$$
\cT_{\cA}(f^{\mathrm{in}},S)(t,x,v)
\leq\|f^{\mathrm{in}}\|_{L^\infty} +\int_0^t\|S(s)\|_{L^\infty}\,ds \leq M_\infty.
$$
\end{proof}

\begin{lemma}
\label{lem:L1}
Let $\mathcal A\subset[0,T]\times\mathbb R^d$ be relatively open and satisfy
$[0,T]\times D_0\subseteq\mathcal A$. Then, for almost every $t\in(0,T)$,
\begin{equation}
\label{eq:L1}
\|\cT_{\cA}(f^{\mathrm{in}},S)(t)\|_{L^1} \leq\|f^{\mathrm{in}}\|_{L^1}+\int_0^t\|S(s)\|_{L^1}\,ds.
\end{equation}
\end{lemma}

\begin{proof}
Extend $f^{\mathrm{in}}$ and $S$ by zero outside $\Omega$ in the spatial variable.
For the initial contribution, nonnegativity and the change of variables $y=x-tv$ give
\begin{align*}
&\int_\Omega\int_{\R^d}
\ind_{\cA}(t,x)\ind_{\{\tau_{\cA}(t,x,v)=0\}}
\ind_{D_0}(x-tv)f^{\mathrm{in}}(x-tv,v)\,dv\,dx
\\
&\qquad\leq \int_{\R^d}\int_{\R^d}f^{\mathrm{in}}(x-tv,v)\,dx\,dv = \|f^{\mathrm{in}}\|_{L^1}.
\end{align*}
For the source contribution, Tonelli's theorem, nonnegativity, and the change of variables
$y=x-(t-s)v$ yield
\begin{align*}
&\int_\Omega\int_{\R^d}\ind_{\cA}(t,x) \int_{\tau_{\cA}(t,x,v)}^t S\bigl(s,x-(t-s)v,v\bigr)\,ds\,dv\,dx
\\
&\qquad\leq
\int_0^t\int_{\R^d}\int_{\R^d} S\bigl(s,x-(t-s)v,v\bigr)\,dx\,dv\,ds = \int_0^t\|S(s)\|_{L^1}\,ds.
\end{align*}
Adding the two estimates proves \eqref{eq:L1}.
\end{proof}

\begin{lemma}
\label{lem:orderT}
Let $\cA_1,\cA_2\subset[0,T]\times\R^d$ be relatively open and suppose that
$[0,T]\times D_0\subseteq\cA_1\subseteq\cA_2$. Then
\begin{equation}
\label{eq:T-order}
\cT_{\cA_1}(f^{\mathrm{in}},S) \leq \cT_{\cA_2}(f^{\mathrm{in}},S) \quad \text{a.e.}
\end{equation}
\end{lemma}

\begin{proof}
By Lemma~\ref{lem:time-order},
$$
\tau_{\cA_2}(t,x,v)\leq\tau_{\cA_1}(t,x,v)
$$
for every $(t,x,v)$.

Fix $(t,x,v)$. If $(t,x)\notin\cA_1$, both contributions to
$\cT_{\cA_1}$ vanish. Suppose instead that $(t,x)\in\cA_1$. Then
$(t,x)\in\cA_2$.
Since $S \geq 0$ and, by Lemma~\ref{lem:time-order},
$\tau_{\cA_2}\leq\tau_{\cA_1}$,
$$
\int_{\tau_{\cA_1}(t,x,v)}^t S\bigl(s,x-(t-s)v,v\bigr)\,ds
\leq \int_{\tau_{\cA_2}(t,x,v)}^t S\bigl(s,x-(t-s)v,v\bigr)\,ds.
$$
Moreover, if $\tau_{\cA_1}(t,x,v)=0$, then
$0\leq\tau_{\cA_2}(t,x,v)\leq\tau_{\cA_1}(t,x,v)=0$, so
$\tau_{\cA_2}(t,x,v)=0$. The remaining factors in the initial contribution are identical.
Formula~\eqref{eq:TA} proves \eqref{eq:T-order}.
\end{proof}

In particular, if $\cD[G_1]\subseteq\cD[G_2]$, then
$$
\cT_{G_1}(f^{\mathrm{in}},S) \leq\cT_{G_2}(f^{\mathrm{in}},S) \quad\text{a.e.}
$$

We conclude this section with some properties of the damage operator, detailed in the following lemmas.

\begin{lemma}
\label{lem:K}
Let
$$
f,f_1,f_2\in L^1\bigl((0,T)\times\Omega\times\R^d\bigr), \qquad f,f_1,f_2\geq0\quad\text{a.e.}
$$
Then $\cK f$, $\cK f_1$, and $\cK f_2$ are finite on
$[0,T]\times\overline\Omega$. If $f_1\leq f_2$ almost everywhere, then
$0\leq\cK f_1\leq\cK f_2$ pointwise. Moreover, for every $t\in[0,T]$,
\begin{equation}
\label{eq:Kbound}
\sup_{x\in\overline\Omega}(\cK f)(t,x) \leq\|K\|_{L^\infty}\int_0^t\|f(s)\|_{L^1}\,ds,
\end{equation}
and, for every $x_1,x_2\in\overline\Omega$,
\begin{equation}
\label{eq:Kmodulus}
|(\cK f)(t,x_1)-(\cK f)(t,x_2)| \leq\omega_K(|x_1-x_2|)\int_0^t\|f(s)\|_{L^1}\,ds.
\end{equation}
\end{lemma}

\begin{proof}
The order property follows from the pointwise nonnegativity of $K$, the almost-everywhere
inequality $f_1\leq f_2$, and monotonicity of the Lebesgue integral.

By Lemma~\ref{lem:kernel-representative}, there exists a measurable set
$N_K\subset\Omega\times\R^d$ of measure zero such that
$0\leq K(x,x',v)\leq\|K\|_{L^\infty}$ for every $x\in\overline\Omega$ and every
$(x',v)\notin N_K$. Hence, for every $t\in[0,T]$ and $x\in\overline\Omega$,
$$
0\leq(\cK f)(t,x) \leq\|K\|_{L^\infty}\int_0^t\|f(s)\|_{L^1}\,ds<\infty.
$$
The same argument applies to $f_1$ and $f_2$. Taking the supremum over $x$ proves
\eqref{eq:Kbound}.

Similarly, Assumption~\ref{H5} gives
\begin{align*}
|(\cK f)(t,x_1)-(\cK f)(t,x_2)|
&\leq\int_0^t\int_\Omega\int_{\R^d} |K(x_1,x',v)-K(x_2,x',v)|f(s,x',v)\,dv\,dx'\,ds
\\
&\leq\omega_K(|x_1-x_2|)\int_0^t\|f(s)\|_{L^1}\,ds,
\end{align*}
which proves \eqref{eq:Kmodulus}.
\end{proof}

\begin{lemma}
\label{lem:damage-continuity}
Suppose that $f\geq0$ and
$$
f\in L^1\bigl((0,T)\times\Omega\times\R^d\bigr).
$$
Then the cumulative damage field $F=\cK f$ belongs to
$C([0,T]\times\overline\Omega)$. For every $x\in\overline\Omega$, the map
$t\mapsto F(t,x)$ is nondecreasing.

If, in addition, there exists $M_1\geq0$ such that
\begin{equation}
\label{eq:generic-L1-bound}
\|f(t)\|_{L^1}\leq M_1 \quad\text{for a.e. }t\in(0,T),
\end{equation}
then
\begin{equation}
\label{eq:space-time-modulus}
|F(t,x)-F(s,y)|\leq \|K\|_{L^\infty}M_1|t-s|+TM_1\omega_K(|x-y|)
\end{equation}
for every $s,t\in[0,T]$ and $x,y\in\overline\Omega$.
\end{lemma}

\begin{proof}
The finiteness of $F$ follows from \eqref{eq:Kbound}. Since $Kf\geq0$, for
$0\leq s\leq t\leq T$ and $x\in\overline\Omega$,
$$
F(t,x)-F(s,x) =\int_s^t\int_\Omega\int_{\R^d} K(x,x',v)f(r,x',v)\,dv\,dx'\,dr\geq0.
$$
Thus $F(\cdot,x)$ is nondecreasing. For arbitrary $s,t\in[0,T]$,
$$
|F(t,x)-F(s,x)| \leq\|K\|_{L^\infty} \int_{\min\{s,t\}}^{\max\{s,t\}}\|f(r)\|_{L^1}\,dr.
$$
The absolute continuity of the Lebesgue integral shows that the right-hand side tends to zero
as $|t-s|\to0$, uniformly in $x$.

By \eqref{eq:Kmodulus},
$$
|F(t,x)-F(t,y)| \leq\omega_K(|x-y|)\int_0^T\|f(r)\|_{L^1}\,dr,
$$
and the right-hand side tends to zero as $|x-y|\to0$, uniformly in $t$.
Combining the two estimates gives the joint continuity of $F$ on
$[0,T]\times\overline\Omega$.

Under \eqref{eq:generic-L1-bound}, the preceding estimates yield
\begin{align*}
|F(t,x)-F(s,y)| &\leq |F(t,x)-F(s,x)|+|F(s,x)-F(s,y)|
\\
&\leq \|K\|_{L^\infty}M_1|t-s| + T M_1\omega_K(|x-y|),
\end{align*}
which proves \eqref{eq:space-time-modulus}.
\end{proof}

\section{Global minimal mild solutions}
\label{sec:existence}

The purpose of this section is to construct a global minimal mild solution by monotone iteration.
Starting from the undamaged field $F^0=0$, we alternately solve the absorbing transport problem in the phase generated by the current damage field and update the damage through $\cK$.

\subsection{The iterative procedure}

Set $F^0:=0$. At each step, $f^{k+1}$ denotes a fixed nonnegative measurable representative of the equivalence class
$\cT_{F^k}(f^{\mathrm{in}},S)$. If two such representatives agree almost everywhere,
then each is almost everywhere bounded above by the other. Hence, by Lemma~\ref{lem:K}, their images under $\cK$ coincide pointwise.

For every $k\geq0$, define recursively
$$
f^{k+1}:=\cT_{F^k}(f^{\mathrm{in}},S), \qquad F^{k+1}:=\cK f^{k+1}.
$$
We also write
$$
D_t^k:=D_t[F^k], \qquad \cD^k:=\cD[F^k].
$$
By definition, $[0,T]\times D_0\subseteq\cD^k$ for every $k$.

\begin{lemma}
\label{lem:iteration-well-defined}
For every $k\geq0$, the function $f^{k+1}$ is measurable, nonnegative, and belongs to
$L^1\bigl((0,T)\times\Omega\times\R^d\bigr) \cap L^\infty\bigl((0,T)\times\Omega\times\R^d\bigr)$.

For every $k\geq1$, the cumulative damage field $F^k:=\mathcal K f^k$ belongs to
$C([0,T]\times\overline\Omega)$, and $t\mapsto F^k(t,x)$ is nondecreasing for every
$x\in\overline\Omega$.
\end{lemma}

\begin{proof}
The field $F^0=0$ is continuous. Hence Lemma~\ref{lem:transport-measurable} gives a measurable representative of
$f^1=\cT_{F^0}(f^{\mathrm{in}},S)$. Formula~\eqref{eq:TA} and nonnegativity of the data give $f^1\geq0$.
Lemmas~\ref{lem:Linf} and~\ref{lem:L1} imply that $f^1\in L^\infty$ and
$f^1\in L^\infty(0,T;L^1)$; in particular, $f^1$ is integrable in space--time.
Lemma~\ref{lem:damage-continuity} then implies that $F^1=\cK f^1$ is continuous and nondecreasing in time at every spatial point.

The same argument applies inductively: if $F^k$ is continuous, then Lemmas~\ref{lem:transport-measurable}, \ref{lem:Linf}, and \ref{lem:L1} give the required properties of $f^{k+1}$, and Lemma~\ref{lem:damage-continuity} gives the corresponding properties of $F^{k+1}$.
\end{proof}

\begin{lemma}
\label{lem:iterorder}
For every $k\geq0$,
$$
F^k\leq F^{k+1}, \qquad
\cD^k\subseteq\cD^{k+1}, \qquad f^{k+1}\leq f^{k+2},
$$
where the inequalities for the representatives of the damage field hold pointwise and those for the densities hold almost everywhere.
\end{lemma}

\begin{proof}
Since $F^0=0$ and $f^1\geq0$, Lemma~\ref{lem:K} gives $F^0\leq F^1$.
Thus $\cD^0\subseteq\cD^1$, and Lemma~\ref{lem:orderT} yields
$$
f^1=\cT_{F^0}(f^{\mathrm{in}},S) \leq \cT_{F^1}(f^{\mathrm{in}},S)=f^2
$$
almost everywhere.

Suppose that $F^{k-1}\leq F^k$. Then $\cD^{k-1}\subseteq\cD^k$, and
Lemma~\ref{lem:orderT} gives $f^k\leq f^{k+1}$ almost everywhere.
Lemma~\ref{lem:K} then gives $F^k\leq F^{k+1}$ pointwise, and consequently $\cD^k\subseteq\cD^{k+1}$. The conclusion follows by induction.
\end{proof}

\begin{proposition}
\label{prop:uniform}
For every $k\geq1$,
\begin{equation}
\label{eq:fk-Linfty}
0\leq f^k\leq M_\infty\quad\text{a.e.}
\end{equation}
For almost every $t\in(0,T)$,
\begin{equation}
\label{eq:fk-L1}
\|f^k(t)\|_{L^1} \leq \|f^{\mathrm{in}}\|_{L^1} +t\|S\|_{L^\infty(0,T;L^1)}.
\end{equation}
Furthermore, for every $k\geq1$ and every
$(t,x)\in[0,T]\times\overline\Omega$,
\begin{equation}
\label{eq:Funiform}
0\leq F^k(t,x) \leq \|K\|_{L^\infty}
\left(t\|f^{\mathrm{in}}\|_{L^1} + \frac{t^2}{2}\|S\|_{L^\infty(0,T;L^1)}\right).
\end{equation}
\end{proposition}

\begin{proof}
Estimate~\eqref{eq:fk-Linfty} follows from Lemma~\ref{lem:Linf}.
Lemma~\ref{lem:L1} gives
$$
\|f^k(t)\|_{L^1} \leq\|f^{\mathrm{in}}\|_{L^1}+\int_0^t\|S(s)\|_{L^1}\,ds,
$$
and \eqref{eq:fk-L1} follows from the $L^\infty(0,T;L^1)$ bound on $S$.
Using \eqref{eq:Kbound} and \eqref{eq:fk-L1}, we obtain
\begin{align*}
F^k(t,x) &\leq\|K\|_{L^\infty}\int_0^t\|f^k(s)\|_{L^1}\,ds
\\
&\leq\|K\|_{L^\infty} \int_0^t\left(\|f^{\mathrm{in}}\|_{L^1} + s\|S\|_{L^\infty(0,T;L^1)}\right)\,ds,
\end{align*}
which proves \eqref{eq:Funiform}.
\end{proof}

\subsection{Convergence of the transport operators}

Numerical convergence of the entrance times to zero does not by itself imply convergence of the indicators of the event that the entrance time is exactly zero. The initial contribution therefore requires a stabilization argument.

\begin{lemma}
\label{lem:initial-stabilization}
Let $\{\cA^k\}_{k\geq0}$ be an increasing sequence of relatively open subsets of
$[0,T]\times\R^d$ such that $[0,T]\times D_0\subseteq\cA^k$ for every $k$.
Set $\cA:=\bigcup_{k\geq0}\cA^k$, $\tau_k:=\tau_{\cA^k}$, and
$\tau:=\tau_{\cA}$. Fix $(t,x,v)$ and suppose that
$$
\tau(t,x,v)=0, \qquad x-tv\in D_0.
$$
Then there exists $k_0$ such that
$$
\tau_k(t,x,v)=0 \qquad \text{for every } k\geq k_0.
$$
\end{lemma}

\begin{proof}
The case $t=0$ is immediate. Assume $t>0$ and set $y:=x-tv=X_{t,x,v}(0)$.
Since $y\in D_0$ and $D_0$ is open, continuity of the characteristic gives
$\varepsilon\in(0,t]$ such that
$$
X_{t,x,v}(r)\in D_0 \qquad\text{for every }r\in[0,\varepsilon].
$$
Hence the characteristic portion over $[0,\varepsilon]$ belongs to every $\cA^k$.

Since $\tau(t,x,v)=0$, one has
$$
(r,X_{t,x,v}(r))\in\cA \qquad \text{for every } r\in(0,t].
$$
Therefore the compact set
$$
\mathcal C_\varepsilon := \bigl\{(r,X_{t,x,v}(r)):r\in[\varepsilon,t]\bigr\}
$$
is contained in $\cA=\bigcup_{k\geq0}\cA^k$. Compactness and monotonicity of the cover yield an index $k_0$ such that
$\mathcal C_\varepsilon\subseteq\cA^{k_0}$. Thus, for every $k\geq k_0$, the whole characteristic segment over $[0,t]$ lies in $\cA^k$, and consequently $\tau_k(t,x,v)=0$.
\end{proof}

Combining the convergence of the backward entrance times with
Lemma~\ref{lem:initial-stabilization}, we obtain the following monotone convergence property.

\begin{lemma}
\label{lem:T-monotone-convergence}
Let $\{G^k\}_{k\geq0}$ and $G$ belong to
$C([0,T]\times\overline\Omega)$. Suppose that
$$
\cD[G^k]\subseteq\cD[G^{k+1}] \quad\text{for every }k,
\qquad \cD[G]=\bigcup_{k\geq0}\cD[G^k].
$$
Then, for almost every $(t,x,v)$, the sequence
$\{\cT_{G^k}(f^{\mathrm{in}},S)(t,x,v)\}$ is nondecreasing and
\begin{equation}
\label{eq:T-limit}
\cT_{G^k}(f^{\mathrm{in}},S)(t,x,v)\to \cT_G(f^{\mathrm{in}},S)(t,x,v).
\end{equation}
\end{lemma}

\begin{proof}
Write $\tau_k:=\tau_{G^k}$ and $\tau:=\tau_G$.
By Lemma~\ref{lem:time-stability}, the sequence $\{ \tau_k \}$ is nonincreasing and $\tau_k\to \tau$ pointwise
as $k \to +\infty$.
The phase indicators are nondecreasing and converge pointwise to
$\ind_{\cD[G]}$.

Let $\widetilde S$ be a nonnegative jointly Borel-measurable representative of $S$.
By Lemma~\ref{lem:representatives}, replacing $S$ by $\widetilde S$ does not change the transport operators almost everywhere. Lemma~\ref{lem:characteristic-null-sets} and Fubini's theorem imply that, for almost every $(t,x,v)$,
$$
0\leq\widetilde S\bigl(s,x-(t-s)v,v\bigr) \leq\|S(s)\|_{L^\infty}
$$
for almost every $s\in(0,t)$. Hence the function
$h(s):=\widetilde S(s,x-(t-s)v,v)$ is integrable on $(0,t)$ for almost every
$(t,x,v)$.

For every such $(t,x,v)$, the indicators
$\ind_{\{\tau_k(t,x,v)<s<t\}}$ increase to
$\ind_{\{\tau(t,x,v)<s<t\}}$ for almost every $s$. The monotone convergence theorem gives
$$
\int_{\tau_k(t,x,v)}^t h(s)\,ds \to \int_{\tau(t,x,v)}^t h(s)\,ds.
$$
Together with convergence of the phase indicators, this proves convergence of the source contributions.

For the initial contribution, if $\tau(t,x,v)>0$, then
$\tau_k(t,x,v)\geq\tau(t,x,v)>0$ for every $k$, so all initial contributions vanish. If $\tau(t,x,v)=0$ but $x-tv\notin D_0$, they vanish because of the factor $\ind_{D_0}(x-tv)$. Finally, if $\tau(t,x,v)=0$ and $x-tv\in D_0$, Lemma~\ref{lem:initial-stabilization}, applied with $\cA^k=\cD[G^k]$, shows that $\tau_k(t,x,v)=0$ for all sufficiently large $k$. Thus the initial contributions eventually coincide with the limiting one.

Lemma~\ref{lem:orderT} gives monotonicity of the complete transport terms, and \eqref{eq:T-limit} follows.
\end{proof}

\subsection{Existence and minimality}

The following theorem is the main result of this section.
\begin{theorem}
\label{thm:existence}
Under Assumptions~\ref{H1}--\ref{H5}, for every finite $T>0$, there exists a global mild solution $(F,f)$ on $[0,T]$ in the sense of Definition~\ref{def:solution}. It satisfies
$$
0\leq f\leq M_\infty\qquad\text{a.e.},
$$
\begin{equation}
\label{eq:solution-L1}
\|f(t)\|_{L^1} \leq \|f^{\mathrm{in}}\|_{L^1} +t \|S\|_{L^\infty(0,T;L^1)}
\quad\text{for a.e. }t\in(0,T),
\end{equation}
and
\begin{equation}
\label{eq:solution-F-bound}
0\leq F(t,x) \leq\|K\|_{L^\infty} \left(t\|f^{\mathrm{in}}\|_{L^1}
+ \frac{t^2}{2}\|S\|_{L^\infty(0,T;L^1)}\right)
\end{equation}
for every $(t,x)\in[0,T]\times\overline\Omega$.
Moreover, $F\in C([0,T]\times\overline\Omega)$, and
$t\mapsto F(t,x)$ is nondecreasing for every $x\in\overline\Omega$.

The solution is minimal among all nonnegative mild solutions: if
$(\widehat F,\widehat f)$ is any other global mild solution, then
\begin{equation}
\label{eq:minimality}
F\leq\widehat F\quad\text{pointwise}, \qquad f\leq\widehat f\quad\text{a.e.}
\end{equation}
\end{theorem}

\begin{proof}
By Lemma~\ref{lem:iterorder}, the sequence $\{f^k\}_{k\geq1}$ is nondecreasing almost everywhere. Since the set of indices is countable, there exists a null set $N_f$ outside which all inequalities
$f^k\leq f^{k+1}$ and all bounds $0\leq f^k\leq M_\infty$ hold simultaneously. Define
$$
f(t,x,v):=\lim_{k\to\infty}f^k(t,x,v)
$$
outside $N_f$, and set $f=0$ on $N_f$. Then $f$ is measurable and
$0\leq f\leq M_\infty$ almost everywhere.

By Fubini's theorem, for almost every $t$ the section $(N_f)_t$ has measure zero.
For every such $t$, monotone convergence gives
$$
\|f(t)\|_{L^1}=\lim_{k\to\infty}\|f^k(t)\|_{L^1}.
$$
Passing to the limit in \eqref{eq:fk-L1} proves \eqref{eq:solution-L1}.
Integration in time gives $f\in L^1$, while the pointwise bound gives $f\in L^\infty$.

Define the cumulative damage field $F:=\cK f$ pointwise as in \eqref{eq:canonical-F}. For fixed
$(t,x)$, the monotone convergence theorem gives
\begin{equation}
\label{eq:Fk-pointwise}
F^k(t,x)\to F(t,x).
\end{equation}
Passing to the limit in \eqref{eq:Funiform} proves \eqref{eq:solution-F-bound}.
Lemma~\ref{lem:damage-continuity} implies that
$F\in C([0,T]\times\overline\Omega)$ and is nondecreasing in time at every spatial point.

Since $F^k$ is nondecreasing and $F^k\to F$ pointwise as $k\to +\infty$, Lemma~\ref{lem:strict-threshold-stability} gives
\begin{equation}
\label{eq:phase-union}
\cD[F]=\bigcup_{k\geq0}\cD[F^k].
\end{equation}
By construction,
$f^{k+1}=\cT_{F^k}(f^{\mathrm{in}},S)$. Lemma~\ref{lem:T-monotone-convergence} and \eqref{eq:phase-union} imply
$$
f^{k+1}\to\cT_F(f^{\mathrm{in}},S)
\quad\text{a.e.}
$$
Since $f^{k+1}\to f$ almost everywhere, we conclude that
$$
f=\cT_F(f^{\mathrm{in}},S) \quad \text{a.e.}
$$
Thus $(F,f)$ is a global mild solution.

Let $(\widehat F,\widehat f)$ be another global mild solution. Since
$\widehat F=\cK\widehat f\geq0=F^0$, one has
$\cD[F^0]\subseteq\cD[\widehat F]$. Lemmas~\ref{lem:orderT} and \ref{lem:K} give
$$
f^1\leq\widehat f\quad\text{a.e.}, \qquad F^1\leq \widehat F \quad \text{pointwise}.
$$
If $F^k\leq\widehat F$, then
$\cD[F^k]\subseteq\cD[\widehat F]$, and the same two lemmas yield
$$
f^{k+1}\leq\widehat f\quad\text{a.e.}, \qquad F^{k+1}\leq\widehat F \quad \text{pointwise}.
$$
Induction and passage to the limit give \eqref{eq:minimality}.
\end{proof}

\begin{corollary}
\label{cor:phase-integrity-monotonicity}
Let $(F,f)$ be any global mild solution.
Then, for every $0\leq s\leq t\leq T$, $D_s[F]\subseteq D_t[F]$.
Moreover,
$u(t,x)\leq u(s,x)$ for every $x\in\overline\Omega$.
\end{corollary}

\begin{proof}
By Lemma~\ref{lem:damage-continuity},
$F(s,x)\leq F(t,x)$ for every $x\in\overline\Omega$.
Therefore,
$$
\{x\in U_0:F(s,x)>1\} \subseteq \{x\in U_0:F(t,x)>1\}.
$$
Since
$$
D_r[F] = D_0\cup\{x\in U_0:F(r,x)>1\},
$$
it follows that $D_s[F]\subseteq D_t[F]$.

Finally, since $r\mapsto[1-r]_+$ is nonincreasing,
$$
[1-F(t,x)]_+ \leq [1-F(s,x)]_+.
$$
Multiplication by $\ind_{U_0}(x)$ yields
$$
u(t,x)\leq u(s,x)
$$
for every $x\in\overline\Omega$.
\end{proof}

\begin{remark}
The construction proves minimality but not uniqueness. In particular, the argument does not exclude the existence of other mild solutions above the minimal one.
\end{remark}

\section{Numerical simulations}\label{s:numerical}
    \subsection{Single-particle erosion law and kernel specification}
The general kernel function $K$ introduced in Section~\ref{sec:model} can be specified to reflect characteristic features of erosion. Here, we focus on the volume of material removed by a single impacting particle. 
In \cite{slikkerveer1998erosion,wensink2002fabrication,fan2009modeling}, it has been shown that the volume removed by a single blasting particle is well expressed in terms of the kinetic energy of the particle and physical properties of the target material, including Young's modulus, Vickers hardness, and fracture toughness. 
In particular, power-law expressions involving the particle kinetic energy have been found to agree well with experimental results \cite{slikkerveer1998erosion}. 
With this in mind, we assume that the volume removed by a single impacting particle is given by
\begin{align}\label{eq:Vp}
    V_{\mathrm{p}} = A\left( \frac{1}{2}m|v|^{2} \right)^{1+\kappa}, 
\end{align}
where $A$ is a given constant depending on the material properties, $m$ is the mass of a single particle, and $v$ is its velocity. 
In the following numerical simulations, we consider the two-dimensional case ($d = 2$), so that $V_{\mathrm{p}}$ is interpreted as an area. 

Our aim is to examine whether the kinetic model can reproduce the erosion rate corresponding to this power-law relation. To this end, we consider the following function:
\begin{align}\label{eq:K-expr}
     K_{\#}(x,x',v) = \mathds{1}_{|x-x'|<r}\tilde{K}_{\#}(v), \quad
    \tilde{K}_{\#}(v) = cA\left( \frac{1}{2}m|v|^{2} \right)^{1+\kappa}\frac{1}{m}\frac{|v|}{r}\frac{1}{r^{2}}, 
\end{align}
where $r$ is the interaction radius and $c$ is a dimensionless constant. 
We have accounted for the following to reach the form of $\tilde{K}_{\#}$ in \eqref{eq:K-expr}. 
First, the decrease in the integral of the residual integrity
\begin{equation}\label{eq:scrU-def}
\mathscr{U}(t):=\int_{U_0} u(t,x)\,dx
\end{equation}
over time can be regarded as the time evolution of the globally removed volume, which, neglecting the truncation in the definition of $u$ in \eqref{eq:uF},
is approximated by
\begin{align*}
\mathscr U(0)-\mathscr U(t)\approx
\int_{U_0}\int_{0}^{t}\int_{D_{s}}\int_{\mathbb{R}^{2}} & K_{\#}(x,x',v)f(s,x',v)\,dv\,dx'\,ds\,dx \\
= & 
\int_{U_0}\int_{0}^{t}\int_{D_{s}}\int_{\mathbb{R}^{2}} \mathds{1}_{|x-x'|<r}\tilde{K}_{\#}(v)f(s,x',v)\,dv\,dx'\,ds\,dx.\notag
\end{align*}
Considering an impingement of a single particle with speed $|v|$, due to $\mathds{1}_{|x-x'|<r}$ and $f$ in the integrand, the integrations over the phase space ($v$ and $x'$), time ($s$), and position ($x$) approximately amount to $m$, $r/|v|$, and $r^{2}$, respectively. 
Multiplying $\tilde{K}_{\#}(v)$ by these factors and the proportionality constant $c^{-1}$, and equating the result to the volume $V_{\mathrm{p}}$, we obtain the expression for $\tilde{K}_{\#}(v)$ in \eqref{eq:K-expr}. 

Although the kernel $K_{\#}$ in \eqref{eq:K-expr} is physically motivated by the single-particle erosion law \eqref{eq:Vp}, it does not satisfy Assumption~\ref{H4} in Section~\ref{sec:transport}, since $\tilde{K}_{\#}\notin L^\infty(\R^d)$. 
Moreover, the presence of the indicator function $\mathds{1}_{|x-x'|<r}$ prevents $K_{\#}$ from satisfying Assumption~\ref{H5}. 
Strictly speaking, therefore, $K_{\#}$ falls outside the class of kernels considered in Sections~\ref{sec:model}--\ref{sec:existence}. 
However, a kernel satisfying these assumptions can be constructed by a slight modification of $K_{\#}$. 
To see this, consider
\begin{align}\label{eq:K-reg-expr}
     K(x,x',v) = \bigg[ \mathds{1}_{|x-x'|<r} \ast \left (\frac{1}{2\pi\epsilon^{2}}\right)\exp\bigg(-\frac{|x-x'|^{2}}{2\epsilon^{2}}\bigg) \bigg]\tilde{K}_{\#}(v)\chi_{V_{m},V_{M}}(|v|). 
\end{align}
Here, $*$ denotes convolution with respect to $x$; $\epsilon$, $V_m$, and $V_{M}$ are given constants such that $\epsilon > 0$ and $0<V_m<V_M$; and $\chi_{V_m,V_M}(|v|)$ is a smooth function taking values in $[0,1]$ such that
$$
\chi_{V_m,V_M}(|v|)=1 \text{ for }|v| \le V_m \text{ and }
\chi_{V_m,V_M}(|v|)=0 
\text{ for }|v| \ge V_M.
$$
By choosing $V_m$ larger than the particle speeds considered in the simulations, the prescribed velocity dependence of $\tilde{K}_{\#}$ is preserved over the relevant range, while $K$ becomes bounded and thus satisfies Assumption~\ref{H4}.
Furthermore, the convolution with the Gaussian ensures that $K$ satisfies Assumption~\ref{H5}.
For sufficiently small $\epsilon$, this regularized kernel approximates $K_{\#}$ well.
In the numerical simulations presented below, we use the simpler kernel $K_{\#}$ for computational convenience.
Since the numerical scheme employs finite time steps and spatial cells, and evaluates the spatial contribution of the kernel through cell averages, the effect of the difference between the two kernels on the numerical results is expected to be negligible for sufficiently small $\epsilon$.
This provides a justification for using $K_{\#}$ as a computational approximation to a kernel satisfying Assumptions~\ref{H4} and \ref{H5}.

\subsection{Numerical method}
We use a particle method to numerically solve the kinetic model with the kernel function $K_{\#}$ given by \eqref{eq:K-expr}.
Denoting by $W$ the number of real particles represented by one simulation particle, the velocity distribution function is approximated as 
\begin{align}\label{eq:VDF-app}
    f(t,x,v) \approx \sum_{k = 1}^{N_{t}} Wm\delta(x-x_{k}(t))\delta(v-v_{k}(t)), 
\end{align}
where $N_{t}$ is the number of simulation particles at time $t$, and $x_{k}(t)$ and $v_{k}(t)$ are the position and velocity of the $k$-th simulation particle, respectively. 
Note that $N_{t}$ simulation particles represent $N_{t}W$ real particles.

The domain $U_{0}$ initially occupied by the substrate is discretized into small square cells $C_{\ell}$ of side length $\Delta x$, where $\ell$ denotes the cell index. 
Then, the value of the cumulative damage $F$ at the center $x_{\mathrm{c}}^{(\ell)}$ of each cell $C_{\ell}$ is approximated as its mean value over the cell: 
\begin{equation}\label{eq:u-app}
    F(t,x_{\mathrm{c}}^{(\ell)}) \approx \frac{1}{|C_{\ell}|}\int_{C_{\ell}}F(t,x)\,dx. 
\end{equation}
From \eqref{eq:canonical-F}, \eqref{eq:VDF-app}, and \eqref{eq:u-app}, we have the following update rule for $F$ after integrating with respect to $v$ and $x'$: 
\begin{align}
\label{eq:u-update}
&F(t^{(n+1)},x_{\mathrm{c}}^{(\ell)}) \! = \! F(t^{(n)},x_{\mathrm{c}}^{(\ell)})+\frac{Wm}{|C_{\ell}|} \int_{t^{(n)}}^{t^{(n+1)}}\! \sum_{k = 1}^{N_{s}}\bigg(\int_{C_{\ell}}\!\mathds{1}_{|x-x_{k}(s)|<r}\,dx\bigg)\tilde{K}_{\#}(v_{k}(s))\,ds, 
\end{align}
where $t^{(n)} = n\Delta t$ and $\Delta t$ denotes the time step. 
In \eqref{eq:u-update}, $\int_{C_{\ell}}\mathds{1}_{|x-x_{k}(s)|<r}\,dx$ is equal to the area of intersection of the cell $C_{\ell}$ and the disk of radius $r$ with its center at $x_{k}(s)$. 

Now, the process of the computation is summarized as follows: 
\begin{enumerate}
\item Advect each particle over a half time step $\Delta t/2$: $x_{k}^{*} = x_{k}(t^{(n)}) + v_{k}(t^{(n)})(\Delta t/2)$. \label{item:first-process}
\item According to \eqref{eq:u-update} and using the midpoint particle positions $\{x_{k}^{*}\}$ and velocities $v_{k}(t^{(n)})$, 
monitor $1-F(t^{(n+1)},x_{\mathrm{c}}^{(\ell)})$, and compute $u(t^{(n+1)},x_{\mathrm{c}}^{(\ell)})$ and $\mathscr{U}(t^{(n+1)})$ based on \eqref{eq:uF} and \eqref{eq:scrU-def}. 
\item Advect each particle over the remaining half time step:\par\noindent $x_{k}(t^{(n+1)}) = x_{k}^{*} + v_{k}(t^{(n)})(\Delta t/2)$.
\item Remove any particle whose new position $x_{k}(t^{(n+1)})$ lies outside the simulation domain $\Omega$ or within a cell with 
$F(t^{(n+1)},x_{\mathrm{c}}) < 1$.
\item Generate simulation particles in the domain $\Omega$ according to the source \par\noindent $\Delta t S(t^{(n+1)},x,v) \,dx\,dv$. \label{item:last-process}
\end{enumerate}
We repeat the processes (\ref{item:first-process})--(\ref{item:last-process}) for $n = 0, 1, \dots $. 
Note that the area of intersection in \eqref{eq:u-update} is approximately computed by dividing $C_{\ell}$ into smaller squares, each with a side length equal to $2^{-M}\Delta x$, and summing the areas of the sub-squares contained within the disk. 
If the disk boundary passes through a divided small square, its area is counted as half of its full value.
Setting $M = 3$ is enough to achieve an accuracy within the resolution of the figures presented below.

\subsection{Single beam impinging normally on a substrate\label{sec:single-beam}}
\subsubsection{Setting and some notations\label{sec:setting-single-beam}}
Consider a square substrate with side length $a$. 
Initially, there are no particles.
Real particles are generated at a rate of $N_{\mathrm{g}}$ per unit time, distributed uniformly at random over a square with side length $a_{\mathrm{s}}$, located at a distance $h$ above the substrate, and impinge normally on the substrate at a given fixed speed $v_{\mathrm{n}}$. 
We investigate the erosion of the substrate based on the kinetic model with $K_{\#}$ given by \eqref{eq:K-expr}. 

Below, we explicitly specify the domains $\Omega$ and $U_{0}$ and the source function $S$ for the above setting. 
Denoting $x = (x_{1}, x_{2})$, we take the coordinates so that the upper surface of the substrate at the initial time is tangent to the line $x_{2} = 0$ and the line $x_{1} = 0$ passes through the centers of the substrate and the source region.
Then, we take the simulation domain $\Omega$ to be a rectangle with width $L_{1}$ and height $L_{2}$, namely,
$$
\Omega = \{ (x_{1},x_{2}) : -L_{1}/2 < x_{1} < L_{1}/2,\ -L_{2}+h_{\Omega} < x_{2} < h_{\Omega} \},
$$
where $h_{\Omega}$ is the vertical distance between the substrate and $\Omega$. 
The substrate at the initial time is expressed as $$U_{0} = \left \{ (x_{1},x_{2}) : -a/2 \le x_{1} \le a/2,\ -a \le x_{2} \le 0 \right \}.$$ 
The initial distribution $f^{\mathrm{in}}$ is zero, whereas the source function $S$ is written as $$S(t,x,v) = m(N_{\mathrm{g}}/a_{\mathrm{s}}^{2})\mathds{1}_{-a_{\mathrm{s}}/2 < x_{1} < a_{\mathrm{s}}/2,\ h < x_{2} < h + a_{\mathrm{s}}}\delta(v-v_{\mathrm{s}})$$ with $v_{\mathrm{s}} = (0,-v_{\mathrm{n}})$. 
Strictly speaking, this source function $S$ does not satisfy Assumption~\ref{H3} because of the presence of the Dirac delta function. 
However, replacing it with a narrow Gaussian distribution,
$(2\pi\sigma^2)^{-1}\exp(-|v-v_{\mathrm{s}}|^2/2\sigma^2)$,
with $\sigma>0$, yields a source function satisfying Assumption~\ref{H3}.
For sufficiently small $\sigma$, the effect of this replacement on the numerical results is expected to be negligible at the resolution considered here.
We therefore retain the original expression for $S$ in the numerical simulations.
The same consideration applies to the source function used in Section~\ref{sec:setting-two-beams}.

For the convenience of specifying the parameters in the simulation and presenting the results, let us introduce some reference quantities in advance.
We denote by $v_{\mathrm{ref}}$ a reference speed. 
Then, consider an area etched by $N_{\mathrm{ref}}$ real particles with speed $v_{\mathrm{ref}}$, and take its square root as the reference length $L_{\mathrm{ref}}$: $$L_{\mathrm{ref}} = \sqrt{N_{\mathrm{ref}}A(mv_{\mathrm{ref}}^{2}/2)^{1+\kappa}}.$$ 
The reference time is taken to be $t_{\mathrm{ref}} = L_{\mathrm{ref}}/v_{\mathrm{ref}}$. 

We have set the substrate size to be $a = 40L_{\mathrm{ref}}$, the number of real particles generated per unit time to be $N_{\mathrm{g}} = 2N_{\mathrm{ref}}/t_{\mathrm{ref}}$, the source region size to be $a_{\mathrm{s}} = 8L_{\mathrm{ref}}$, and the separation distance to be $h = 20L_{\mathrm{ref}}$. 
The constant $\kappa$ in \eqref{eq:K-expr} is taken to be either $0$ or $0.25$. 
These values are chosen since the exponents corresponding to $\kappa$ used in \cite{slikkerveer1998erosion} are $0.23$, $0.247$, or $1/6$, which fall between the chosen values. 
The interaction radius $r$ is fixed at $r = 0.5L_{\mathrm{ref}}$. 
The dimensionless constant $c$ is set to $1.5$, with which we confirmed that $\int_{U_0} u(t,x)\,dx$ decreases by nearly $1.0 L_{\mathrm{ref}}^{2}$ when $N_{\mathrm{ref}}/W$ simulation particles with identical initial position and speed $v_{\mathrm{ref}}$ are injected. 
The particle speed is varied as $v_{\mathrm{n}}/v_{\mathrm{ref}} = 0.125$, $0.25$, or $0.5$. 
The dimensions of the domain $\Omega$ are specified by $L_{1} = 80L_{\mathrm{ref}}$, $L_{2} = 90L_{\mathrm{ref}}$, and $h_{\Omega} = 30L_{\mathrm{ref}}$. 

\subsubsection{Results\label{sec:result-single-beam}}
We show below the results obtained using $N_{\mathrm{ref}}/W = 1000$ (which determines the number of simulation particles), the cell size $\Delta x/L_{\mathrm{ref}} = 0.05$, and the time steps $\Delta t/t_{\mathrm{ref}} = 0.005$, $0.0025$, or $0.00125$ respectively for the cases of $v_{\mathrm{n}}/v_{\mathrm{ref}} = 0.125$, $0.25$, or $0.5$. 
For the results shown below, numerical convergence has been judged within the error that
does not affect the statements in this section through the comparison with the results obtained with other parameters corresponding to coarser space and time discretizations and fewer simulation particles. 

\begin{figure}[bt]
    \centering 
\includegraphics[width=0.85\textwidth]{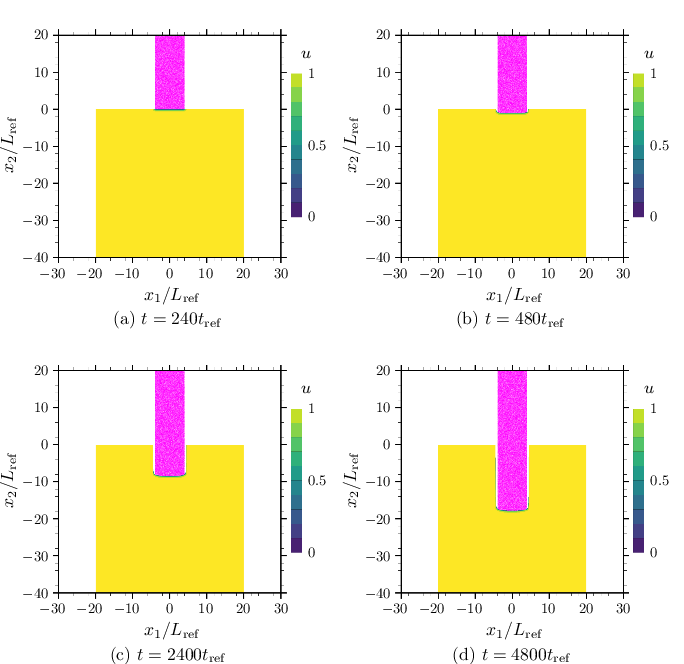}
\caption{
Colormap of $u(t,x)$ for $\kappa = 0$ and $v_{\mathrm{n}} = 0.125v_{\mathrm{ref}}$. Particle positions are shown by markers as well (appearing as a thick straight line extending from the upper part of the region). 
(a)~$t = 240t_{\mathrm{ref}}$, 
(b)~$t = 480t_{\mathrm{ref}}$, 
(c)~$t = 2400t_{\mathrm{ref}}$, and 
(d)~$t = 4800t_{\mathrm{ref}}$. 
The undamaged part of the substrate ($u=1$) is shown by yellow, while the color bar shows the range $0 \le u < 1$.
}
\label{fig:u_single-beam-vn-0.125}
\end{figure}

\begin{figure}[bt]
    \centering 
\includegraphics[width=0.85\textwidth]{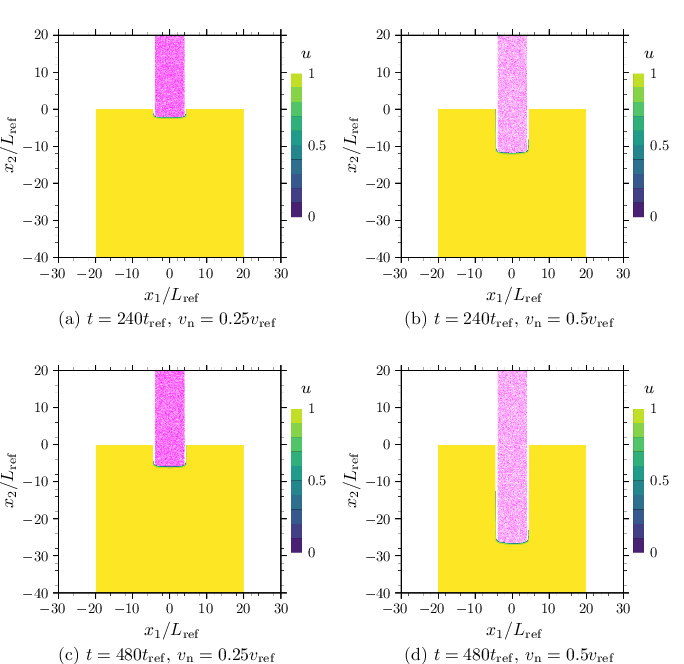}
\caption{
Colormap of $u(t,x)$ for $\kappa = 0$. 
(a)~$t = 240t_{\mathrm{ref}}$, $v_{\mathrm{n}} = 0.25v_{\mathrm{ref}}$, 
(b)~$t = 240t_{\mathrm{ref}}$, $v_{\mathrm{n}} = 0.5v_{\mathrm{ref}}$, 
(c)~$t = 480t_{\mathrm{ref}}$, $v_{\mathrm{n}} = 0.25v_{\mathrm{ref}}$, and 
(d)~$t = 480t_{\mathrm{ref}}$, $v_{\mathrm{n}} = 0.5v_{\mathrm{ref}}$. 
}
\label{fig:u_single-beam-vn-faster}
\end{figure}

Figure~\ref{fig:u_single-beam-vn-0.125} shows the profiles of the residual integrity $u(t,x)$ for $\kappa = 0$ and $v_{\mathrm{n}} = 0.125v_{\mathrm{ref}}$ at several times $t$. 
The positions of simulation particles are also shown.
The substrate initially undergoes partial weakening ($0<u<1$) without the formation of an eroded region accessible to particles [see panel (a)]. 
Indeed, the set of markers indicating the particle positions terminates at the surface $x_{2} = 0$. 
The void region accessible to particles grows as time goes on. 
Figure~\ref{fig:u_single-beam-vn-faster} shows the corresponding profiles for $\kappa = 0$ at higher incident particle speeds, $v_{\mathrm{n}} = 0.25v_{\mathrm{ref}}$ and $0.5v_{\mathrm{ref}}$.
As expected, the substrate is eroded faster for larger $v_{\mathrm{n}}$, as can be seen by comparing Figs.~\ref{fig:u_single-beam-vn-0.125} and \ref{fig:u_single-beam-vn-faster}. 
These results illustrate how the kinetic model describes the accumulation of damage and the subsequent evolution of the eroded region.

\begin{figure}[bt]
    \centering 
\includegraphics[width=0.55\textwidth]{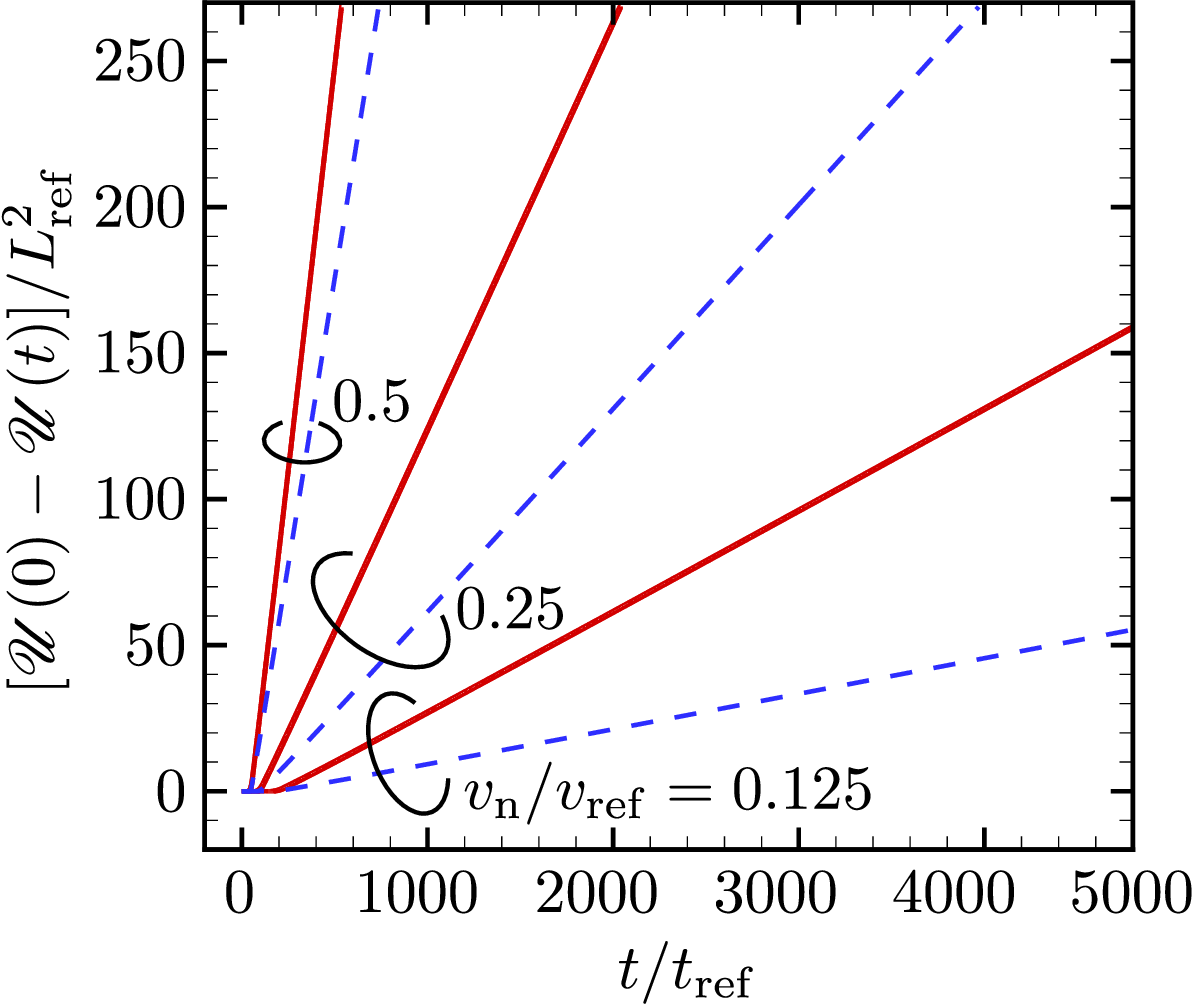}
\caption{
$t/t_{\mathrm{ref}}$ vs 
$[\mathscr U(0) - \mathscr U(t)]/
L_{\mathrm{ref}}^{2}$. 
The results for $\kappa = 0$ and $0.25$ are indicated by the solid and dashed lines, respectively. 
}
\label{fig:115_plot_Udx_vs_t}
\end{figure}

Next, we examine the time evolution of the total loss of residual integrity, defined by 
$[\mathscr{U}(0)-\mathscr U(t)]$ over time both for $\kappa = 0$ and $0.25$ in Fig.~\ref{fig:115_plot_Udx_vs_t}. 
For the particle speeds considered here, the larger exponent corresponding to $\kappa=0.25$ results in slower erosion than in the case of $\kappa=0$.
After an initial transient period (e.g., around $t/t_{\mathrm{ref}} < 200$ for $v_{\mathrm{n}}/v_{\mathrm{ref}} = 0.125$ where the plotted curves exhibit slight curvature), the plotted curves appear nearly linear in time. 
To confirm this in an objective manner, we carry out linear regression on each data set. 
For each pair of $\kappa$ and $v_{\mathrm{n}}$, we choose an appropriate upper bound of the time interval $t_{\mathrm{f}}$ and duration $\Delta$, and fit a line $a^{(m)} (t/t_{\mathrm{ref}}) + b^{(m)}$ to the numerical data on the interval $[t_{\mathrm{f}}-m\Delta, t_{\mathrm{f}}]$ ($m = 1, 2, 3$). 
For example, for $\kappa = 0$ and $v_{\mathrm{n}} = 0.25v_{\mathrm{ref}}$, taking $t_{\mathrm{f}} = 1800t_{\mathrm{ref}}$ and $\Delta = 450t_{\mathrm{ref}}$, we have $a^{(3)} = 0.1390$, $b^{(3)} = -15.04$, $a^{(2)} = 0.1392$, $b^{(2)} = -15.38$, $a^{(1)} = 0.1392$, and $b^{(1)} = -15.36$. 

\begin{figure}[bt]
    \centering 
    \includegraphics[width=0.55\textwidth]{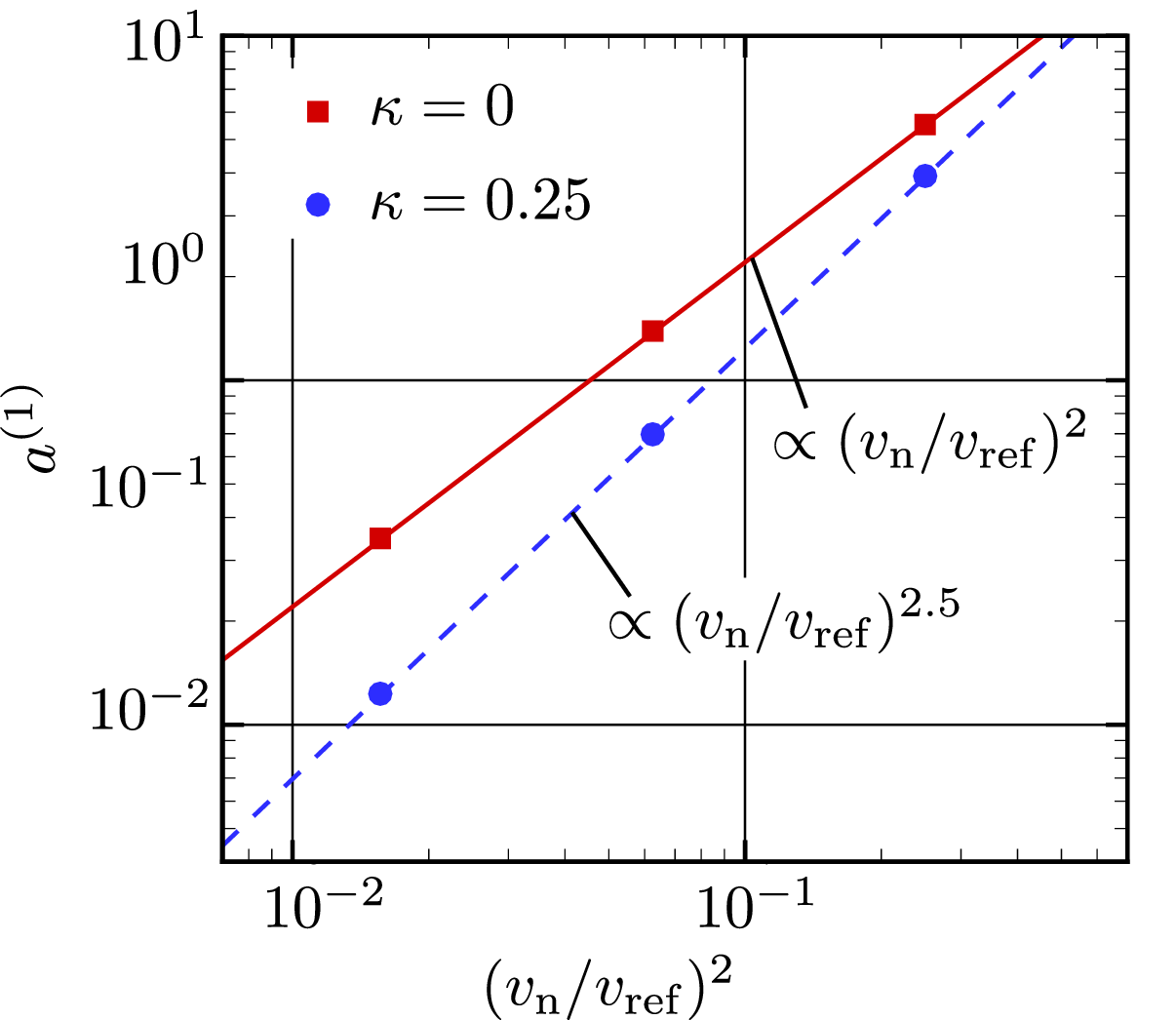}
\caption{
$(v_{\mathrm{n}}/v_{\mathrm{ref}})^{2}$ vs the slope $a^{(1)}$ read from the graphs in Fig.~\ref{fig:115_plot_Udx_vs_t} by the least-squares method. 
}
\label{fig:201_slopes_G3}
\end{figure}

Figure~\ref{fig:201_slopes_G3} shows the slopes obtained from the linear regression for different values of the incident particle speed. 
This implies that the amount of erosion over a given time interval is proportional to $v_{\mathrm{n}}^{2(1+\kappa)}$, which matches the dependence on the particle kinetic energy assumed in the single-particle erosion law \eqref{eq:Vp}. 
Note that the particle mass flux is the same across all speeds $v_{\mathrm{n}}$, since $N_{\mathrm{g}}$ is kept fixed. 
Thus, the erosion rate proportional to the powers of kinetic energy can be reproduced by the present kinetic model, demonstrating the role of the kernel function $K_{\#}$ or $K$ in determining the erosion behavior.


\begin{figure}[bt]
    \centering 
\includegraphics[width=0.85\textwidth]{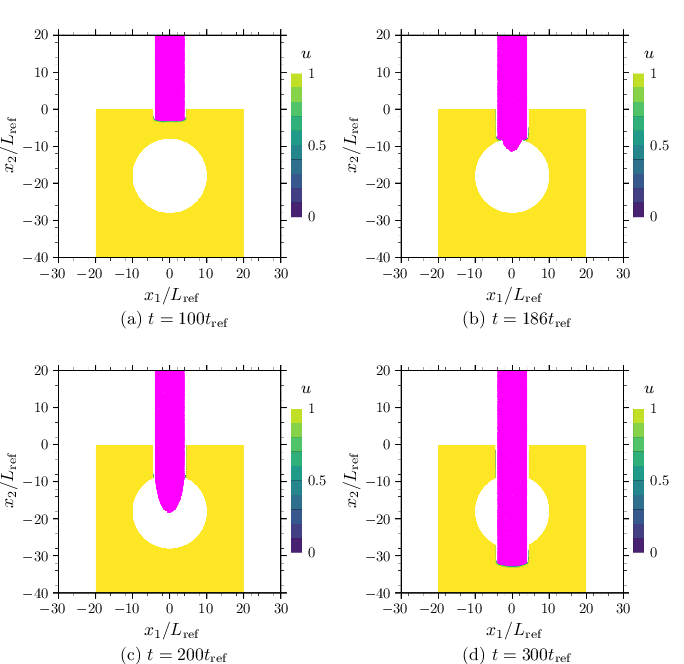}
\caption{
Colormap of $u(t,x)$ for a substrate containing a circular hole, with $\kappa = 0$ and $v_{\mathrm{n}} = 0.5v_{\mathrm{ref}}$. Particle positions are shown by markers as well (appearing as a thick straight line extending from the upper part of the region). 
(a)~$t = 100t_{\mathrm{ref}}$, 
(b)~$t = 186t_{\mathrm{ref}}$, 
(c)~$t = 200t_{\mathrm{ref}}$, and 
(d)~$t = 300t_{\mathrm{ref}}$.  
}
\label{fig:u_single-beam-with-a-hole}
\end{figure}

Finally, we consider a substrate containing a circular hole. 
The setting is the same as before, except that the substrate initially contains a circular hole centered at $(0,-18L_{\mathrm{ref}})$ with radius $10L_{\mathrm{ref}}$. 
Figure~\ref{fig:u_single-beam-with-a-hole} shows the profiles of the residual integrity $u(t,x)$, together with the positions of simulation particles, for $\kappa=0$ and $v_{\mathrm{n}}=0.5v_{\mathrm{ref}}$. 
Initially, the substrate is eroded from its upper surface, as in the case without a hole. 
As erosion proceeds, the eroded region reaches the hole, and the two initially disconnected particle-accessible regions become connected. 
The particles then enter the hole and subsequently erode the substrate beneath it. 
This example illustrates how the present kinetic model can describe erosion processes involving changes in the topology of the particle-accessible region.

\subsection{Two beams impinging obliquely on a substrate\label{sec:setting-two-beams}}
We consider the same setting as in Section~\ref{sec:setting-single-beam}, except that particles are injected obliquely toward the substrate from two square source regions, one on the left and the other on the right.
Both beams are directed toward the central region of the substrate, each making an angle of $45^\circ$ with the vertical direction. 
The source function $S$ is given by $$S(t,x,v) = m(N_{\mathrm{g}}/a_{\mathrm{s}}^{2})\mathds{1}_{-a_{\mathrm{s}}/2-x_{\mathrm{s}} < x_{1} < -x_{\mathrm{s}}+a_{\mathrm{s}}/2,\ h < x_{2} < h + a_{\mathrm{s}}}\delta(v-v_{\mathrm{s}}^{\mathrm{L}})$$
$$\qquad\qquad\qquad+ m(N_{\mathrm{g}}/a_{\mathrm{s}}^{2})\mathds{1}_{-a_{\mathrm{s}}/2+x_{\mathrm{s}} < x_{1} < x_{\mathrm{s}}+a_{\mathrm{s}}/2,\ h < x_{2} < h + a_{\mathrm{s}}}\delta(v-v_{\mathrm{s}}^{\mathrm{R}})$$ with $v_{\mathrm{s}}^{\mathrm{L}} = (v_{\mathrm{b}}/\sqrt{2},-v_{\mathrm{b}}/\sqrt{2})$ and $v_{\mathrm{s}}^{\mathrm{R}} = (-v_{\mathrm{b}}/\sqrt{2},-v_{\mathrm{b}}/\sqrt{2})$.
Here, $v_{\mathrm{b}}$ denotes the speed of the incident particles and $2x_{\mathrm{s}}$ is the distance between the centers of the two squares. 

We have set the substrate size to be $a = 40L_{\mathrm{ref}}$, the number of real particles generated per unit time in each square to be $N_{\mathrm{g}} = 0.2N_{\mathrm{ref}}/t_{\mathrm{ref}}$, the source region size to be $a_{\mathrm{s}} = 3L_{\mathrm{ref}}$, and the separation distances to be $h = 20L_{\mathrm{ref}}$ and $x_{\mathrm{s}} = 24L_{\mathrm{ref}}$. 
The parameters of the kernel function are set to $\kappa=0$, $r=0.5L_{\mathrm{ref}}$, and $c=1.5$. 
The particle speed is fixed at $v_{\mathrm{b}}=0.5v_{\mathrm{ref}}$. The simulation domain is specified by $L_1=80L_{\mathrm{ref}}$, $L_2=90L_{\mathrm{ref}}$, and $h_{\Omega}=30L_{\mathrm{ref}}$. 
We show below the results obtained using $N_{\mathrm{ref}}/W = 1000$ (which determines the number of simulation particles), the cell size $\Delta x/L_{\mathrm{ref}} = 0.05$, and the time step $\Delta t/t_{\mathrm{ref}} = 0.00125$. 

\begin{figure}[bt]
    \centering 
    \includegraphics[width=0.85\textwidth]{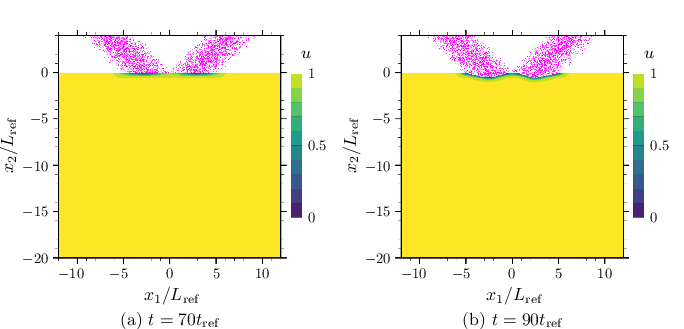}
    \caption{
Colormap of $u(t,x)$ in the problem described in Section~\ref{sec:setting-two-beams}. Particle positions are shown by markers as well. 
(a)~$t = 70t_{\mathrm{ref}}$ and (b)~$t = 90t_{\mathrm{ref}}$. 
}
\label{fig:u_two-beams-former}
\end{figure}

\begin{figure}[bt]
    \centering 
    \includegraphics[width=0.85\textwidth]{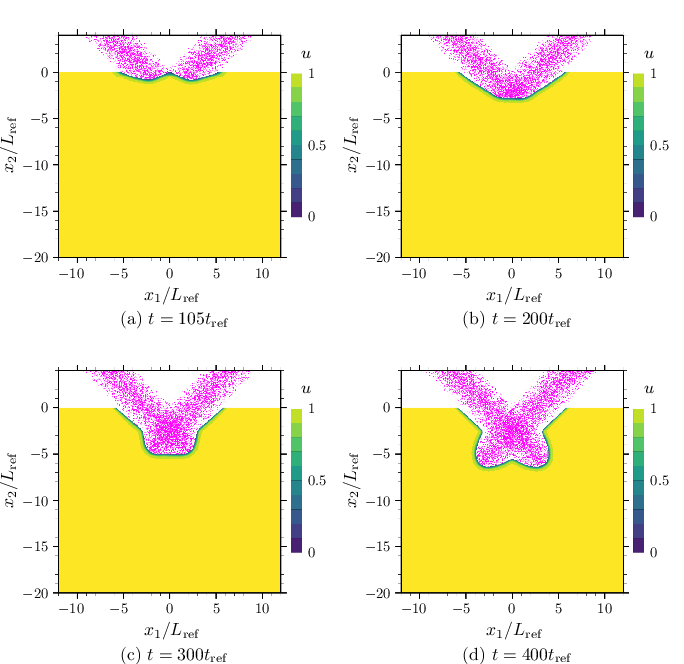}
    \caption{
Colormap of $u(t,x)$. 
(a)~$t = 105t_{\mathrm{ref}}$, 
(b)~$t = 200t_{\mathrm{ref}}$, 
(c)~$t = 300t_{\mathrm{ref}}$, and
(d)~$t = 400t_{\mathrm{ref}}$. 
See also the caption of Fig.~\ref{fig:u_two-beams-former}. 
}
\label{fig:u_two-beams-latter}
\end{figure}

Figures~\ref{fig:u_two-beams-former} and \ref{fig:u_two-beams-latter} show the time evolution of the spatial distribution of the residual integrity $u(t,x)$, together with the positions of the simulation particles. 
In the early stage, the two particle beams erode the substrate independently, forming two separate depressions on its surface (Fig.~\ref{fig:u_two-beams-former}). Toward the end of this stage, the portion of the substrate between the two depressions develops a sharp tip [Fig.~\ref{fig:u_two-beams-latter}(a)]. As erosion proceeds further, the two depressions eventually merge into a single eroded region, which continues to extend into the substrate [Figs.~\ref{fig:u_two-beams-latter}(b)--(d)]. 

An interesting feature of this process is that the two beams interact indirectly through the evolution of the region accessible to the particles. 
Although there are no direct interactions between particles, erosion caused by one beam changes the region accessible to particles from the other beam, thereby affecting the subsequent erosion process. 
The present example illustrates how the kinetic model can describe such an interaction through the evolution of the substrate geometry.


\section{Conclusion}
\label{sec:discussion}

In this paper, we have introduced a kinetic weakening--erosion model in which cumulative damage, threshold-driven accessibility, and absorption along particle
characteristics are treated separately.
A topological formulation of the backward entrance time and a monotone iteration scheme yield a global minimal
mild solution on every prescribed finite time interval.

The main advantage of this approach is its robustness with respect to nonsmooth interfaces and topology changes.
The construction, indeed, does not require a classical interface normal, a boundary parametrization, or a
prescribed interface velocity, and remains meaningful when the geometry of the accessible phase is irregular.

However, several questions remain open because the present construction relies essentially on positivity, monotonicity, and phase inclusion.
These properties provide existence and minimality, but they do not yield two-sided stability estimates or control the behavior at the threshold.
They may also fail when the erosion law depends on the evolving geometry, as in genuinely angle-dependent models.
Consequently, uniqueness, compatibility across different time horizons, angle-dependent erosion, and geometric regularity of the threshold set may 
require methods beyond the  framework developed here.

In the numerical simulations presented in this paper, we have demonstrated that the model can reproduce a power-law dependence of the erosion rate on the kinetic energy of impacting particles, as one of the basic characteristics of erosion. 
However, the kernel function employed in these simulations is not designed to capture other features of more realistic erosion processes. 
Examples include the dependence of the erosion rate on the particle impact angle and the transition from brittle to ductile erosion that may occur when the impacting particles are very small \cite{sheldon1966ductile,wensink2002fabrication}. 
Constructing kernel functions that incorporate these effects and examining them through numerical simulations remain topics for future research.


\medskip

\noindent
\textit{Acknowledgements.} The authors acknowledge support from the JSPS-MAEDI bilateral joint research projects (SAKURA programs: Nos. JPJSBP120173202, JPJSBP120193219, and JPJSBP120263205) for the initiation and continuation of the collaboration that led to this study. The specific blasting model presented in this paper was conceptualized during the second author's stay at Kyoto University, supported by a JSPS International Fellowship for Research in Japan (S24088). The second author also thanks GNFM-INdAM for its support and Kyoto University for its warm hospitality during his stay. The first and third authors also thank ESILV for its hospitality during their stay under the support of the SAKURA programs.

\medskip

\noindent
\textit{Conflicts of interest.}
The authors state that there is no conflict of interest.

\medskip

\noindent
\textit{Declaration of generative AI and AI-assisted technologies in the writing process.}
During the preparation of this work, the authors used Gemini (Google) to assist with correcting English language usage, and improving text clarity and organization. 
After using this tool, the authors reviewed and edited the content as needed and take full responsibility for the content of the publication.

\medskip

\noindent
\textit{Data availability.}
The data that support the findings of this article are not publicly available. 
The data are available from the authors upon reasonable request.

\addcontentsline{toc}{section}{Bibliography}
\bibliography{biblio}

\end{document}